\documentclass[halfline]{IET}
\usepackage{cite}
\usepackage{graphicx}
\usepackage{epstopdf}
\usepackage{color,soul}
\graphicspath{{Figures/}}
\begin{document}

\title{Stable Galerkin Finite Element Scheme for the Simulation of Problems Involving Conductors Moving Rectilinearly in Magnetic Fields}

\author[1,*]{Sethupathy Subramanian}

\author[2]{Udaya Kumar}
\affil[1,2]{ Department of  Electrical Engineering, Indian Institute of Science, Bangalore 560012, India}
\affil[*]{sethupathys87@gmail.com}

\abstract{For the simulation of rectilinearly moving conductors across a magnetic field, the Galer-kin finite element method (GFEM) is generally employed. The inherent instability of GFEM is very often addressed by employing Streamline upwinding/Petrov-Galerkin (SU/PG) scheme. However, the SU/PG solution is known to suffer from distortion at the boundary transverse to  the velocity and the remedial measures suggested in fluid dynamics literature are computationally demanding. Therefore, simple alternative schemes are essential. In an earlier effort, instead of conventional finite-difference based approach, the numerical instability was analyzed using the Z-transform. By employing the concept of pole-zero cancellation, stability of the numerical solution was achieved by a simple restatement of the input magnetic flux in terms of associated vector potential. This approach, however, is restricted for input fields, which  vary only along the direction of the velocity.
To overcome this, the present work proposes a novel approach 
in which the input field is restated as a weighted elemental average.
The stability of the proposed scheme is proven analytically for both 1D and  2D cases. The error bound for the small oscillations remnant at intermittent Peclet numbers is also deduced. Using suitable numerical simulations, all the theoretical deductions are verified.}

\maketitle

``This paper is a postprint of a paper submitted to and accepted for publication in \textit{IET Science, Measurement \& Technology} and is subject to Institution of Engineering and Technology Copyright. The copy of record is available at IET Digital Library''

\newpage
\section{Introduction}
Accurate evaluation of the induced currents and the resulting magnetic fields is very essential in the design of equipment, working on electromagnetic induction principle.
Among such equipment, this work basically concerns with the ones involving a rectilinear movement of conducting media under applied magnetic field.
The ready examples are electromagnetic flow meters, linear induction motors and eddy current brakes. The governing equations for steady state electromagnetic fields in such devices are \cite{reducedb} \cite{fmbase}, 
\begin{equation} \label{eqge2}
\nabla \cdot (\sigma \nabla \phi) - \nabla \cdot (\sigma ~ {\bf{u}} \times \nabla \times {\bf{A}}) =  \nabla \cdot (\sigma ~ \bf{u} \times \bf{B_{a}})
\end{equation}	
\begin{equation} \label{eqge1}
\sigma \nabla \phi ~-~  (\nabla \cdot \dfrac{1}{\mu} \nabla) {\bf{A}} - \sigma~ {\bf{u}} \times \nabla \times {\bf{A}} = \sigma~ {\bf{u}} \times {\bf{B_{a}}}
\end{equation}	
 where,  $\bf A$ is the vector potential associated with reaction magnetic field $\bf b_r$, $\phi$ is the scalar potential field arising out of current flow, ${\bf u}$ is the velocity of the moving conductor, $\sigma$ is the electrical conductivity,  $\mu$  is the magnetic permeability and $\bf{B_a}$ is the applied magnetic field. The source of the applied field is considered to be unaffected by the reaction field.
 
 The analytical solution of the governing equations (\ref{eqge2}) and (\ref{eqge1}) is rather difficult and hence numerical techniques are generally employed. For this, the Galerkin finite element method (GFEM), which is a widely employed numerical method across different disciplines, appears to be the best suited option. It is known to give accurate results, when the Peclet number, $Pe = \mu \sigma |{\bf u}| \Delta z /2~<1$ ($\Delta z$ is the element length along the direction of the velocity) \cite{cdbook}\cite{quada1}. However, as the velocity of the conductor/conducting-fluid becomes high, a very fine discretisation involving large number of elements needs to be employed. This would be a practically difficult exercise.
 
 For $Pe > 1$ on the other hand, GFEM is known to suffer from numerical instability along the direction of the velocity. This issue of numerical instability  has been adequately addressed in the fluid dynamics literature \cite{cdbook} mostly for the nodal formulation. Among the methods suggested, the Streamline upwinding/Petrov-Galerkin (SU/PG) scheme \cite{up1}\cite{supg1} is commonly employed in the electromagnetic literature \cite{mc6tf1,mc1lim1,mc2av1,mc2av2,mc3eb1,mc3eb2,mc4ge1, mc5mc1, mc6tf2,mc6tf3}.

It is reported in fluid dynamics literature that, while SU/PG scheme is successful in ensuring stability of the solution, it can lead to localized peaking/distortion  at the boundary transverse to the flow  \cite{discop}. This is a serious issue, which needs careful attention. Numerical experiments carried out in the 2D and 3D version of the flowmeter problem, have clearly confirmed the existence  of such an error at the conductor-air interface transverse to the flow. Remedial measures such as, `discontinuity-capturing' scheme \cite{discop}, Finite Increment Calculus (FIC) \cite{fic2} have  been suggested in fluid dynamics literature, which will overcome this issue.
 However,  these numerical schemes are non-linear and hence demands more computation. In addition, efforts towards applying them for electromagnetic problem is rather scarce.

In view of these, it was deemed necessary to seek within the framework of GFEM, a stable scheme, which is free of errors at the transverse boundary. At this juncture, it will be important to note that the problem under consideration involves rectilinear movement of the conductors. In such problems, at least the conducting region of the problem can be and usually be discretised with graded regular mesh along the direction of the velocity.  In other words, the resulting mesh would be like a stack of layers of different thickness along the direction of the velocity. This aspect will be referred in the later part of the work.
Further for discretisation,  quadrilateral elements for 2D and brick elements for 3D will be considered for the analysis.

Classically, any analysis of the stability is generally carried out with the one-dimensional version of the problem, which is discretised with a regular mesh. For that, the GFEM equation turns out to be difference equation, which is generally solved analytically to investigate on the numerical stability of
the solution. In our earlier work \cite{su1}, this difference equation was translated into Z-transform domain and was analyzed by borrowing tools from the  control system theory. In that, the ratio of the applied magnetic flux density and the reaction magnetic field was defined as the transfer function. A pole of the transfer function was traced to be the source of numerical instability. Then zero in the numerator was brought in for the required pole-zero cancellation, by simply restating the input magnetic field in terms of its vector potential. It ensured absolute stability even at very high velocities.
 This approach however was found to work well, only when the input magnetic field varies along the direction of the velocity.
In practice, as the input magnetic field can have variation orthogonal to  the direction of the velocity, it is necessary to overcome the above limitation and this forms the goal of the present work.

In this paper, firstly a novel stable scheme is devised for the 1D version of the problem. It is then directly applied to the 2D version of the problem and the existence of the stability is analytically shown using the 2D Z-transform analysis. Numerical simulations are carried out to validate the scheme.


\section{Present work}
\subsection{Analysis with the 1D version of the problem}
Following the footsteps of the earlier works \cite{cdbook}, \cite{up1},  \cite{quadp1}, investigation will be carried out first with the 1D version of the problem. The conditions that permit the reduction of the physical problem to its 1D version has been described in \cite{su1} and the corresponding governing equation is,
\begin{equation}\label{eq1d1x}
~~-\dfrac{d^{2}A_{y}}{dz^2} + \mu \sigma u_z \dfrac{dA_y}{dz} = \mu\sigma u_z B_x ~
\end{equation}
where, $A_y$ is the $y$ component of the vector potential, $u_z$ is the velocity of the moving conductor along the $z$ direction and $B_x$ is the input magnetic field. Clearly, (\ref{eq1d1x}) has the same form as the convection-diffusion equation dealt in fluid dynamics \cite{cdbook}, \cite{supg1},\cite{fic1}.  Application of the GFEM  to (\ref{eq1d1x}) leads to difference equation  \cite{cdbook}, \cite{mc6tf2}, which for $n^{th}$ node takes the form,
	\begin{equation}\label{eq1dd2}
	 (-1-Pe)A_{y[n-1]} + 2 A_{y[n]} + (-1+Pe)A_{y[n+1]} = 2Pe\Delta z \Big( \dfrac{B_{x[n-1]}+4B_{x[n]}+B_{x[n+1]}}{6} \Big)
	\end{equation}

It may be noted that, in the early part of \cite{su1} instead of GFEM, the difference approximation was directly employed which results in a slightly different RHS.
	Also in \cite{su1}, it was shown that the use of Z-transform can simplify the analysis of instability.  Following the same, the required analysis will be carried out on the transfer function, defined as the ratio between the vector potential $A_y$ of the reaction magnetic field and the input flux density $B_x$, can be written as,
\begin{equation}
 \dfrac{A_y}{B_x} ~=~ \dfrac{2Pe\Delta z~}{6(-1+Pe)}~~ \dfrac{Z^2+4Z+1}{Z^2~+~\dfrac{2}{-1+Pe}Z~+~\dfrac{-1-Pe}{-1+Pe}} 
\end{equation}

	when $Pe >> 1$
\begin{equation}
 \dfrac{A_y}{B_x} ~\simeq~ \dfrac{\Delta z}{3}~~ \dfrac{(Z+0.27)(Z+3.73)}{(Z~-~1)~(Z~+~1)} 
\end{equation}
As discussed in \cite{su1},	the pole located at -1  is responsible for the oscillation in the computed result \cite{debook1}, \cite{dcbook2}.  It may be recalled that, the above observation is valid for any input magnetic field. Further, $Pe$ forms the true independent variable and not its constituents taken in isolation \cite{cdbook}, \cite{quada1}, \cite{quadp1}. 

The RHS of (\ref{eq1dd2}) can also be viewed as a weighted average of the nodal flux densities, where the averaging is according to the Galerkin formulation. By restating the input magnetic field in terms of vector potential, it was shown in \cite{su1} that necessary zero can be introduced to cancel out the pole. Taking queue from this idea, a consistent modified  weighted average of the input nodal flux densities is sought so as to introduce necessary zero.

The RHS of (\ref{eq1dd2}) with general weighted nodal flux density can be written as, 
\begin{equation} \label{gnrlbavg}
 \tilde{B} = \dfrac{\alpha B_{x[n-1]}+\beta B_{x[n]}+\gamma B_{x[n+1]}}{\alpha+\beta+\gamma} 
\end{equation}
	where $\alpha,~\beta,~\gamma$ are the unknown parameters, which needs to be determined imposing the appropriate constraints. 	
Firstly, total flux into the element must be closely represented. Secondly, the two elements, spanned by the weighing function of the node under consideration, should have equal influence. The last constraint is to obtain the necessary zero (i.e. $Z+1$ term) in the RHS. With these restrictions,  the values of parameters can be evaluated for the linear element as $\alpha=1$, $\beta=2$, $\gamma=1$. Tracing this back to the elemental input, the required averaged input nodal flux densities for the  element $e$ with nodes $n-1$ \& $n$, can be identiﬁed as, 
\begin{equation} \label{ip1d}
B_{e} = \dfrac{B_{x[n-1]}+B_{x[n]}}{2}
\end{equation}

With the above weighted input nodal magnetic field, the equation for the weighing function associated with the $n^{th}$ node reduces to,
\begin{equation}\label{eq1dd3}
	 (-1-Pe)A_{y[n-1]} + 2 A_{y[n]} + (-1+Pe)A_{y[n+1]} = 2Pe\Delta z \Big( \dfrac{B_{x[n-1]}+2B_{x[n]}+B_{x[n+1]}}{4} \Big)
	\end{equation}

The relation  between the reaction field $A_y$ and the input field $B_x$ now takes the form,
	\begin{equation}\label{eq1da1}
		\dfrac{A_y}{B_x} = \dfrac{2Pe\Delta z}{4(-1+Pe)} \dfrac{(Z+1)^2}{(Z-1)\Big(Z-\dfrac{(-1-Pe)}{(-1+Pe)}\Big)}
	\end{equation}
	
	when $Pe >> 1$ the above reduces to

	\begin{equation}\label{eq1da3}
	 \dfrac{A_y}{B_x} ~\simeq~ \dfrac{\Delta z}{2}~~ \dfrac{(Z+1)}{(Z-1)}
	\end{equation}

	\begin{figure}
		\centering
		\subfigure[]{\label{fem1d1:3000} \includegraphics[scale=0.9]{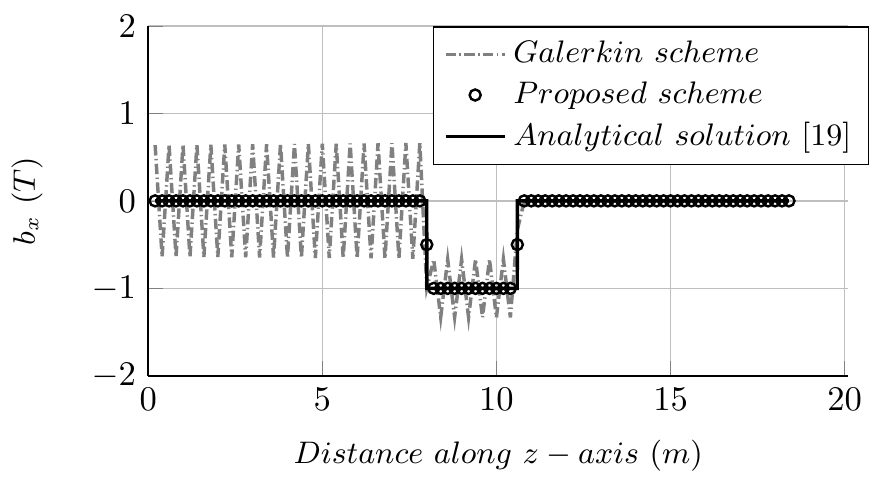}}
		\subfigure[]{\label{fem1d1:0003} \includegraphics[scale=0.9]{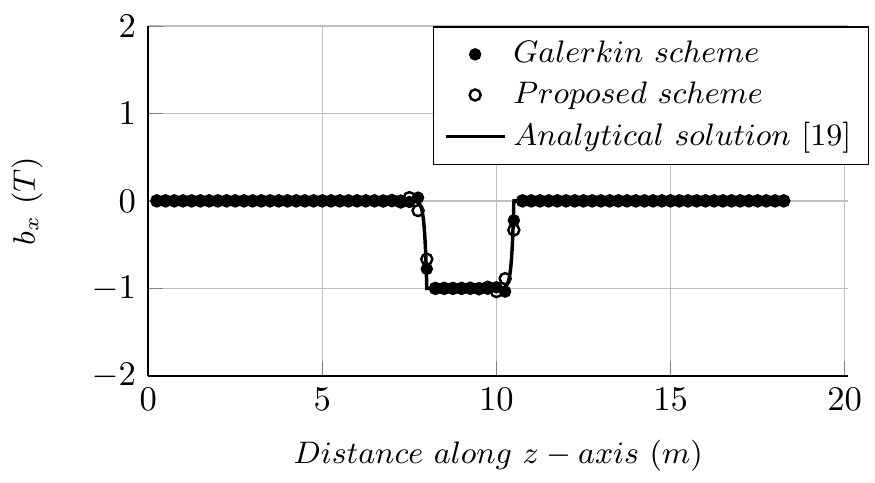}}
		    {\subfigure[]{\label{errplot} 
\includegraphics[scale=0.95]{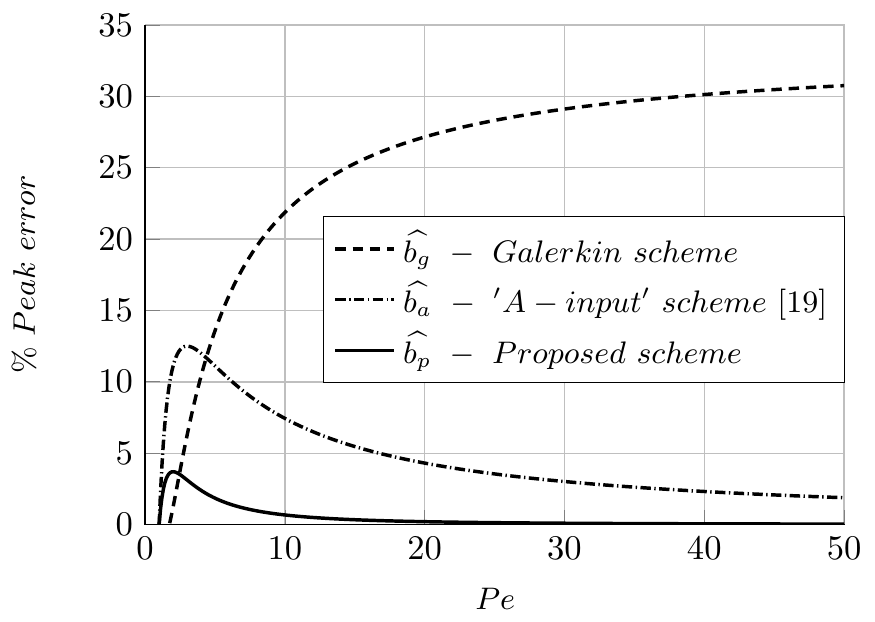}}}
		\caption{1D-FEM solution.  $~$ (a) Pe = 2000, {$\Delta z = 0.2$.}$~~$ (b) Pe = 2, {$\Delta z = 0.25$.} (c) peak error in the numerical solution for a range of  $Pe$}
		\label{fem1d1}
	\end{figure}
It is evident from the above equation that, for $Pe>>1$ the proposed scheme is absolutely stable. Referring to (\ref{eq1da1}) it can be verified that for $Pe<1$, there is no  oscillatory pole. However, for $Pe $ in the range $1$ to $10$ there exists imperfect pole-zero cancellation, leading to small oscillation in the solution. These have been verified with suitable FEM simulation. Sample results are presented in fig. \ref{fem1d1} with  input field ($B_x = B ~when~a \leq z \leq b$) and  boundary conditions $A_y(0) = 0$ and $dA_y/dz | _L = 0$  \cite{su1}. Here, the reaction magnetic field is calculated from  $b_x=-dA_y/dz$.

 In order to quantify the peak amplitude of the small oscillation/error found in the result for the mid-range of Peclet numbers ($1<Pe<10$), analytical solution of the FEM difference-equation is deduced. Details are presented in the Appendix. The error in the simulation result for three different approaches are presented in fig. \ref{errplot}. The peak error in the solution is quantified by (\ref{errp}) and (\ref{errg}) of the Appendix for proposed and Galerkin schemes respectively. 
 
 A maximum error of about 3.7\%  is found to occur at $Pe=2$  in the proposed scheme, which is lower by a factor of 3 with respect to `A-input' scheme  \cite{su1}. This better performance can be attributed to the presence of double `$-1$' zeros in the numerator of (\ref{eq1da1}).

Following the general trend in the pertinent literature on the stability of numerical schemes, required analysis was carried out first with the 1D version of the problem  \cite{cdbook}, \cite{up1},  \cite{quadp1}. However the  main concern in this work is the input magnetic field, which varies even in the direction transverse to the velocity. This calls for suitable analysis with the  2D version of the problem, which will be dealt in the next section.

\subsection{Analysis with the 2D version of the problem}
	\begin{figure}
		\centering
  	    {\subfigure[]{\label{2dprob1} \includegraphics[scale=0.3]{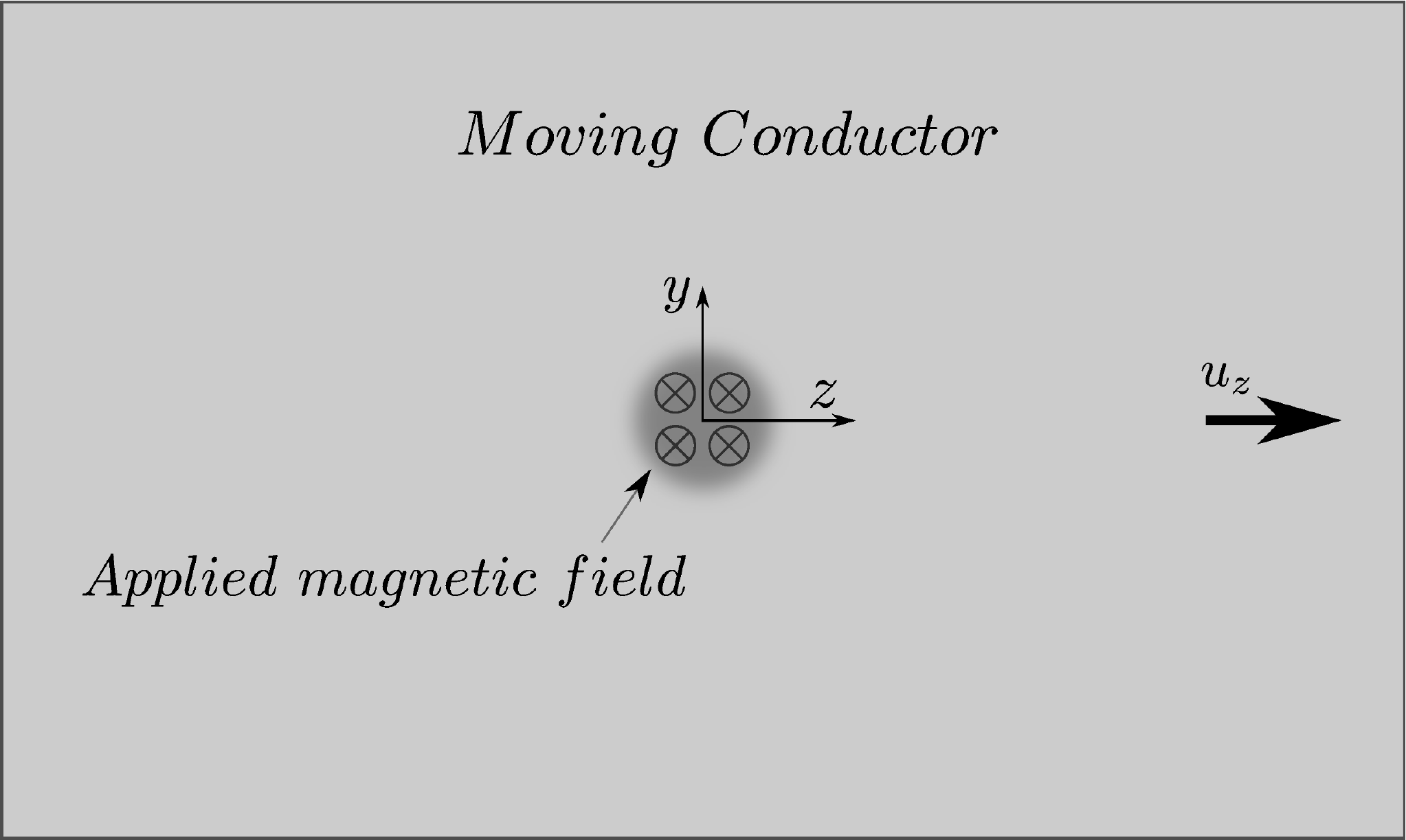}}}
		{\subfigure[]{\label{2ddisc_sq} \includegraphics[scale=0.24]{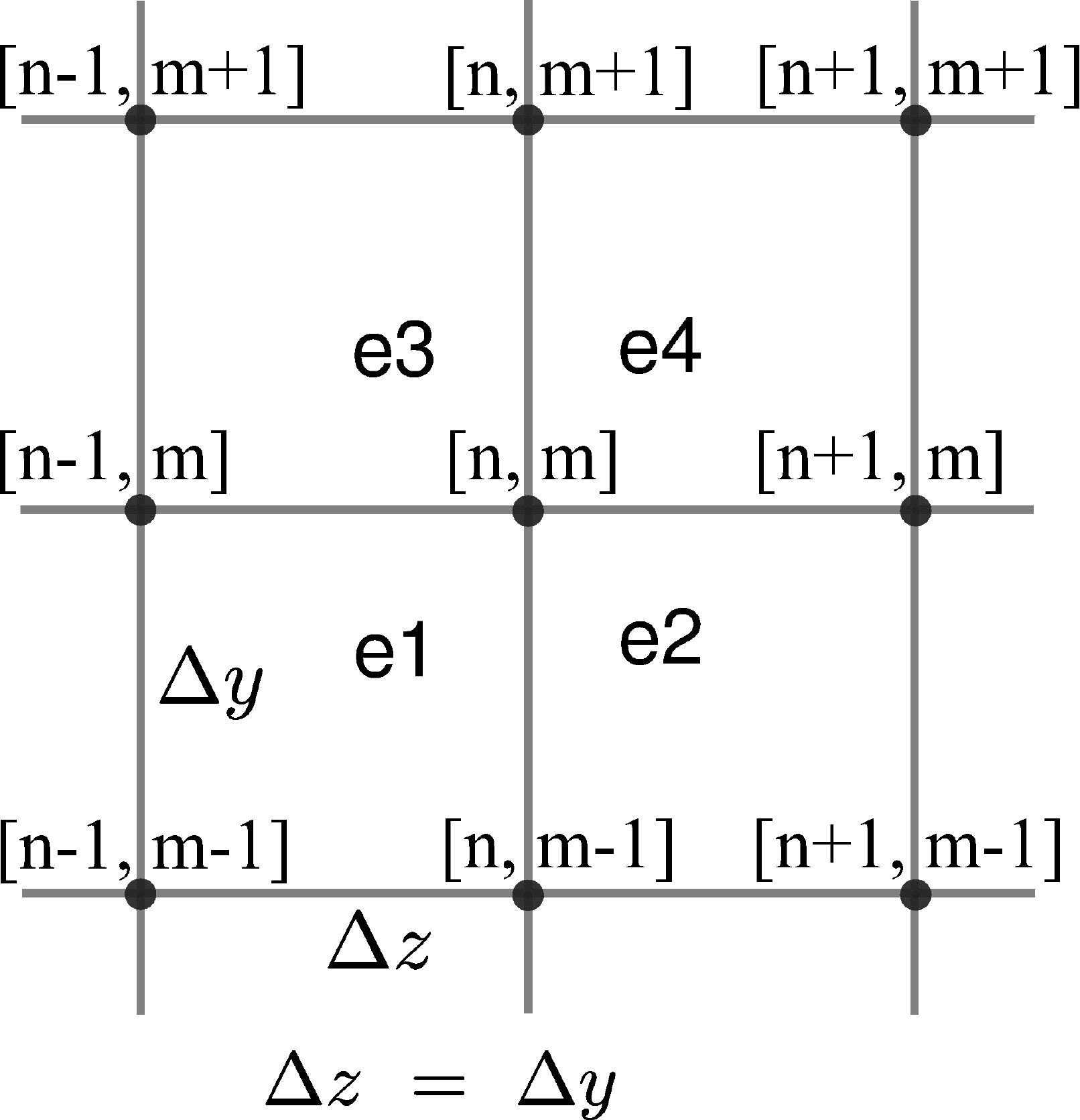}}}
  		{\subfigure[]{\label{2dipsqr} \includegraphics[scale=0.43]{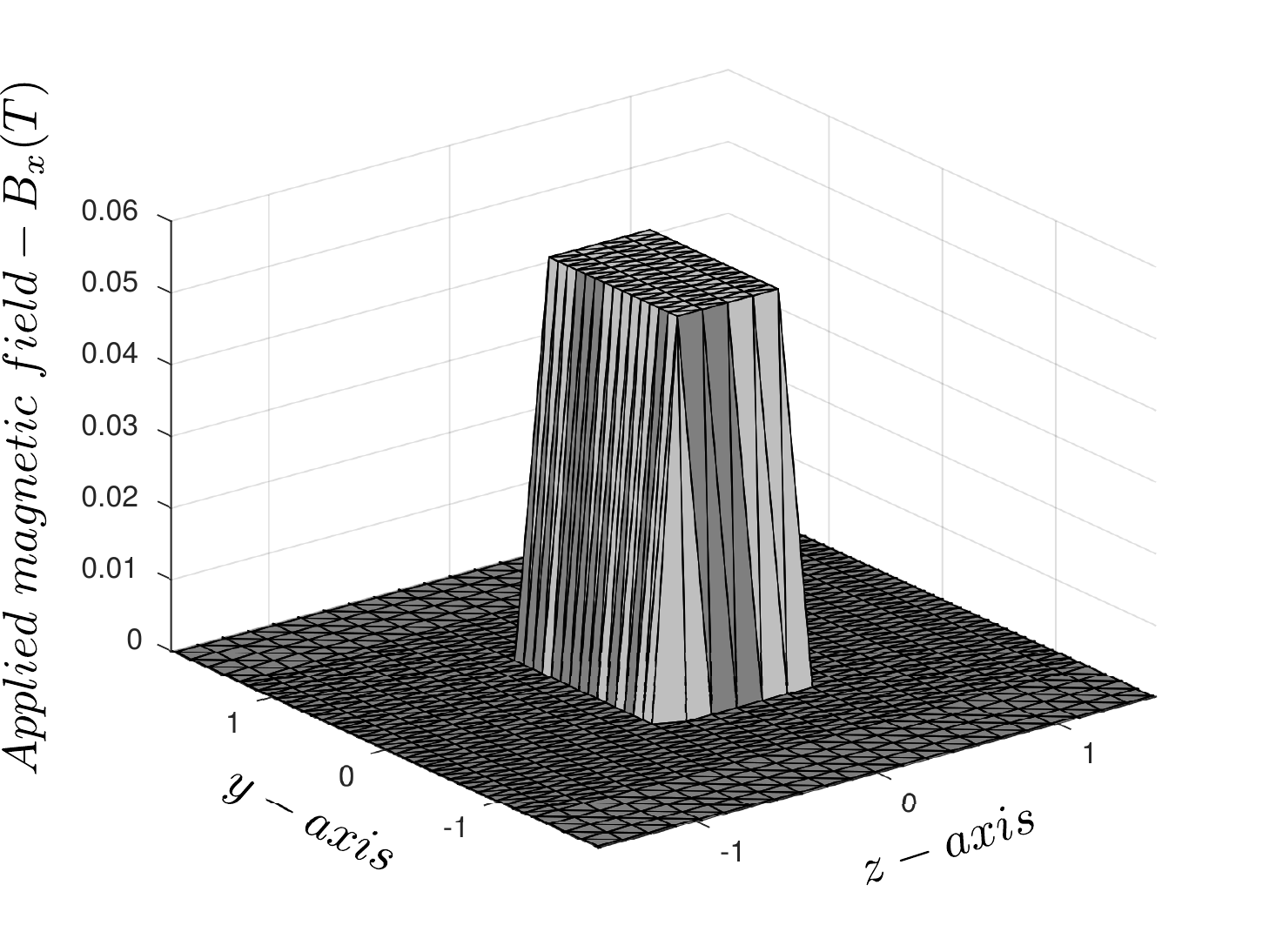}}}
		{\subfigure[]{\label{2dipcrc} \includegraphics[scale=0.43]{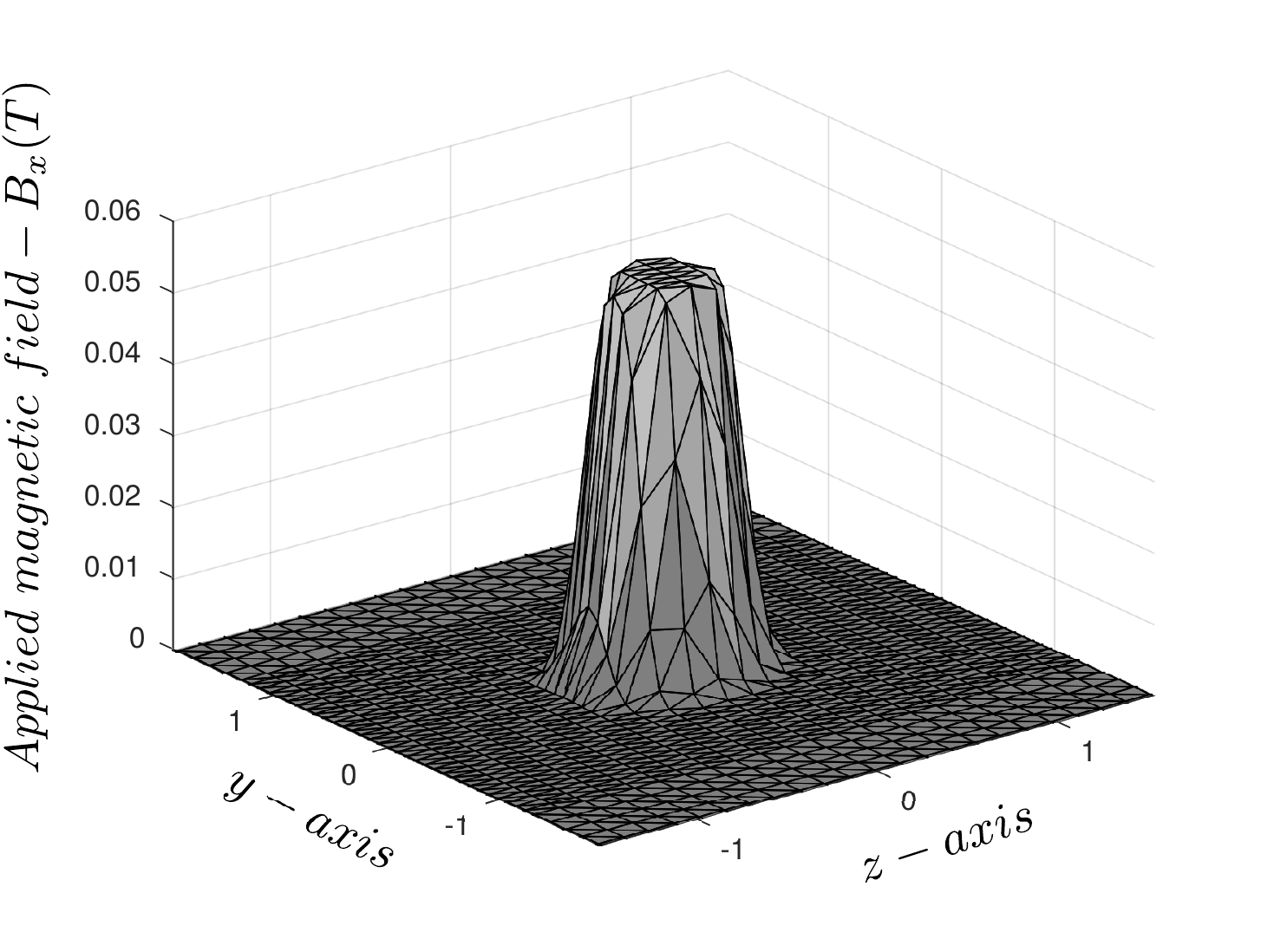}}}
		\caption{2D problem. (a) Schematic of the 2D problem. (b) 2D grid (c) Input magnetic field - rectangular pulse profile. (d) Input magnetic field - smooth circular  profile. }
		\label{fig2dztr}
	\end{figure}
Consider an infinite conducting slab moving along the $z$-axis. A localized, $x$ directed magnetic field, which varies in both $y$ and $z$ direction is applied. Two profiles as shown in Figs. \ref{2dipsqr} and \ref{2dipcrc}  have been considered. For the first one, input field is given by $B_x = B; ~\text{for}~-a \leq z \leq a~\text{and} ~-b \leq y \leq b$ and 
for the second one, $B_x = B; ~\text{for}~ r \leq R$ and $B_x = Be^{-((r-R)/0.5R)^2}~ \text{for} ~r > R$
 with  a smooth Gaussian fall to zero is considered, where $r = \sqrt(y^2+z^2)$.

 The governing equations now take the form,
\begin{equation} \label{ge2dphi}
\dfrac{\partial^{2}\phi}{\partial z^2} + \dfrac{\partial^{2}\phi}{\partial y^2}+u_z \dfrac{\partial^{2}A_y}{\partial y \partial z} - u_z \dfrac{\partial^{2}A_z}{ \partial y^2}= u_z \dfrac{\partial B_x }{\partial y}
\end{equation}
\begin{equation} \label{ge2day}
 \sigma \dfrac{\partial \phi}{\partial y}-\dfrac{1}{\mu}(\dfrac{\partial^{2}A_y}{\partial z^2}+ \dfrac{\partial^{2}A_y}{\partial y^2}) + \sigma u_z \dfrac{\partial A_y}{dz}- \sigma u_z \dfrac{\partial A_z}{dy} = \sigma u_z B_x  
\end{equation}
\begin{equation} \label{ge2daz}
 \sigma \dfrac{\partial \phi}{\partial z}-\dfrac{1}{\mu}(\dfrac{\partial^{2}A_z}{\partial z^2} + \dfrac{\partial^{2}A_z}{\partial y^2}) = 0
\end{equation}

The Z-transform approach will be employed and for this, a regular 2D FEM grid with quadrilateral elements  shown in Fig. \ref{2ddisc_sq} is considered. The corresponding GFEM equation is obtained and for brevity the expanded form is given only for (\ref{ge2daz}),
\begin{equation} \label{2dazdiff}
\begin{split}
\dfrac{Pe}{6 u_z}\big(\phi_{[n+1,m+1]}-\phi_{[n-1,m+1]} +4\phi_{[n+1,m]}-4\phi_{[n-1,m]}+\phi_{[n+1,m-1]}-\phi_{[n-1,m-1]}\big) ~-~~ \\ \dfrac{1}{3}(A_{z[n-1,m-1]}+A_{z[n,m-1]}+A_{z[n+1,m-1]}+A_{z[n-1,m]}-8A_{z[n,m]}\\
+A_{z[n+1,m]}+A_{z[n-1,m+1]}+A_{z[n,m+1]}+A_{z[n+1,m+1]})  = 0
\\
\end{split}
\end{equation}
For this 2D problem, following the literature on multi-dimensional signal processing \cite{2dhb}\cite{2dsi}, Z-transform is applied for both $y$ and $z$ co-ordinates.
Then (\ref{2dazdiff}) reduces to,
\begin{equation} \label{2dazz1}
\begin{split}
\dfrac{Pe}{6 u_z}\big(Z_nZ_m-Z_n^{-1}Z_m +4Z_n-4Z_n^{-1}+Z_nZ_m^{-1}-Z_n^{-1}Z_m^{-1}\big) \phi ~-~~ \\ \dfrac{1}{3}(Z_nZ_m+Z_m+Z_n^{-1}Z_m+Z_n-8+Z_n^{-1}+Z_nZ_m^{-1}+Z_m^{-1}+
Z_n^{-1}Z_m^{-1}) A_z = 0
\\
\end{split}
\end{equation}
where, $Z_n$ and $Z_m$ corresponds to $z$ and $y$ directions respectively. Multiplying (\ref{2dazz1}) by $Z_nZ_m$
\begin{equation} \label{2dazz2}
\begin{split}
\dfrac{Pe}{6u_z}~(Z_n^2Z_m^2-Z_m^2+4Z_n^2Z_m-4Z_m+Z_n^2-1)~\phi~-~ \\ \dfrac{1}{3}(Z_n^2Z_m^2+Z_nZ_m^2+Z_m^2+Z_n^2Z_m-8Z_nZ_m+Z_m+Z_n^2+Z_n+1)~A_z~=~~~0
\end{split}
\end{equation}
which can be written as,
\begin{equation} \label{2dazzc1}
\dfrac{Pe}{6u_z}[Q2]\phi- \dfrac{1}{3}[S1]A_z=0 
\end{equation}
where,
\begin{align}
[S1] &= Z_n^2Z_m^2+Z_nZ_m^2+Z_m^2+Z_n^2Z_m-8Z_nZ_m+Z_m+Z_n^2+Z_n+1 \label{S1} \\
[Q2] &= Z_n^2Z_m^2-Z_m^2+4Z_n^2Z_m-4Z_m+Z_n^2-1 \label{Q2} 
\end{align}

Similarly, the $Z$-transform of the GFEM approximation for (\ref{ge2dphi}) and (\ref{ge2day}) can be reduced to,
\begin{equation} \label{2dphizc1g}
\dfrac{1}{3}[S1]\phi+\dfrac{u_z}{4}[S2]A_y-\dfrac{u_z}{6}[S3]A_z=\dfrac{u_z \Delta z}{12} [Q1]B_x
\end{equation}
\begin{equation} \label{2dayzc1g}
\dfrac{Pe}{6 u_z}[Q1]\phi+\big\{-\dfrac{1}{3}[S1]+\dfrac{Pe}{6}[Q2] \big\}A_y-\dfrac{Pe}{6}[Q1]A_z=\dfrac{Pe \Delta z}{18} [M1]B_x
\end{equation}
where,
\begin{align}
[S2] &=Z_n^2Z_m^2-Z_m^2-Z_n^2+1  \label{S2} \\
[S3] &= Z_n^2Z_m^2+4Z_nZ_m^2+Z_m^2-2Z_n^2Z_m-8Z_nZ_m-2Z_m+Z_n^2+4Z_n+1  \label{S3} \\
[Q1] &= Z_n^2Z_m^2+4Z_nZ_m^2+Z_m^2-Z_n^2-4Z_n-1  \label{Q1}\\
[M1] &= Z_n^2Z_m^2+4Z_nZ_m^2+Z_m^2+4Z_n^2Z_m+16Z_nZ_m+4Z_m+Z_n^2+4Z_n+1  \label{M1} 
\end{align}

The above three, (\ref{2dazzc1}), (\ref{2dphizc1g}), (\ref{2dayzc1g}), describe the nature of the numerical solution of independent field variables ($\phi, A_y, A_z$). For the intended analysis, however, it will be convenient to deal with one equation. It may be recalled that the numerical oscillations are always due to the dominance of the  first order derivative term over the second order derivative \cite{cdbook} and equation for $A_y$ (\ref{ge2day}) possesses this feature.
 This is further confirmed by the numerical experiments, wherein oscillations in the computed $A_y$ is much more dominant than that in the other two variables. In order to obtain equation involving only $A_y$, $\phi$ and $A_z$ are eliminated from (\ref{2dayzc1g})  using equations (\ref{2dazzc1}) and (\ref{2dphizc1g}). For brevity, the intermittent steps are avoided here however will be provided later for the proposed scheme. After simplification for $Pe>>1$, the final equation relating $A_y$ to $B_x$ for the GFEM scheme can be found as,
  \begin{equation} \label{2dgfemtf}
\begin{split}
\dfrac{A_y}{B_x}  &\simeq \dfrac{\Delta z}{3}\dfrac{ (Z_n+0.27)(Z_n+3.7) ~f_1(Z_m)}{(Z_n-1)(Z_n+1)~f_2(Z_m)}\\
\end{split}
\end{equation}
where, $f_1=2(Z_m^2-2Z_m+1)(Z_m^2+4Z_m+1)-3(Z_m^2-1)^2 $; and $f_2 = -(Z_m - 1)^4$. 
 
The numerator and denominator of (\ref{2dgfemtf}) are found to be in \textit{separable form} \cite{2dsi} \cite{2dsb}, 
therefore it is possible to isolate zeros of $Z_n$ and $Z_m$ polynomials. It can be verified that only the polynomial in $Z_n$ appearing in the denominator has a zero at `-1', which is responsible for the numerical oscillation.

For the proposed scheme, dealing with the first order elements, the philosophy adopted earlier for 1D problem can be directly extended to the 2D case. {Recall that the present investigation will be restricted to quadrilateral elements.} 
 Accordingly, in the evaluation of the elemental matrices, the four nodal flux densities are replaced by their arithmetic average. {For illustration, consider the element $e_1$ in Fig. {\ref{2ddisc_sq}}, the elemental averaged flux $B_{e1}$ is given by }
\begin{align}
B_{e1} = \big(B_{x[n-1,m-1]}+B_{x[n,m-1]}+B_{x[n-1,m]}+B_{x[n,m]}\big)/4
\label{2dbe1}
\end{align}
With the above modification, the RHS of ({\ref{2dphizc1g}}) and ({\ref{2dayzc1g}}) will get modified to the form,
\begin{equation} \label{2dphizc1}
\dfrac{1}{3}[S1]\phi+\dfrac{u_z}{4}[S2]A_y-\dfrac{u_z}{6}[S3]A_z=\dfrac{u_z \Delta z}{8} [R1]B_x
\end{equation}
\begin{equation} \label{2dayzc1}
\dfrac{Pe}{6 u_z}[Q1]\phi+\big\{-\dfrac{1}{3}[S1]+\dfrac{Pe}{6}[Q2] \big\}A_y-\dfrac{Pe}{6}[Q1]A_z=\dfrac{Pe \Delta z}{8} [N1]B_x
\end{equation}
where,
\begin{align}
[R1] &= Z_n^2Z_m^2+2Z_nZ_m^2+Z_m^2-Z_n^2-2Z_n-1  \label{R1}  \\
[N1] &= Z_n^2Z_m^2+2Z_nZ_m^2+Z_m^2+2Z_n^2Z_m+4Z_nZ_m+2Z_m+Z_n^2+2Z_n+1   \label{N1} 
\end{align}

 Following the steps described earlier for GFEM, equation for $A_y$ is deduced as follows. Equation (\ref{2dazzc1}) can be rewritten as,
\begin{equation} \label{2dphiz}
\phi=\dfrac{2u_z}{Pe} \dfrac{[S1]}{[Q2]}A_z
\end{equation}
substituting (\ref{2dphiz}) in (\ref{2dphizc1}),
\begin{equation} \label{eq:zt2dbi_phired_f_04}
\begin{split}
\bigg\{ \dfrac{4u_z[S1]^2-u_z Pe [S3][Q2]}{6Pe[Q2]} \bigg \} ~A_z~+~\dfrac{u_z}{4}[S2]~A_y~=~\dfrac{u_z \Delta z}{8} [R1]~B_x
\end{split}
\end{equation}
For $Pe >> 1$, the above can be reduced to,
\begin{equation} \label{2dazz}
A_z \simeq \dfrac{3[S2]}{2[S3]}A_y-\dfrac{3\Delta z[R1]}{4[S3]}~B_x
\end{equation}
For $Pe >>1$, substituting  (\ref{2dphiz}), (\ref{2dazz}) in (\ref{2dayzc1}), the relation between $A_y$ and $B_x$ can be written as,
\begin{equation} \label{2dptf1}
\dfrac{A_{y}}{B_x} \simeq \dfrac{3\Delta z}{2} \dfrac{(~[S3][N1]-[Q1][R1]~)}{(~2[S3][Q2]-3[Q1][S2]~)}
\end{equation}
Expanding the numerator and denominator of (\ref{2dptf1}), we get
\begin{equation} \label{2dptfnm}
[S3][N1]-[Q1][R1] = (Z_n^2+4Z_n+1)(Z_n^2+2Z_n+1)~f_3(Z_m)
\end{equation}
where $f_3(Z_m)=(Z_m^2-2Z_m+1)(Z_m^2+2Z_m+1)-(Z_m^2-1)^2$ and 
\begin{equation} \label{2dptfdn}
2[S3][Q2] - 3[Q1][S2] = (Z_n^2+4Z_n+1)(Z_n^2-1)~f_2(Z_m)
\end{equation}

Similar to the GFEM case, the polynomials in $Z_n$ and $Z_m$ are in separable form and it can be simplified to,
\begin{align*}
\dfrac{A_{y}}{B_x} &\simeq \dfrac{3\Delta z}{2} \dfrac{(Z_n^2+4Z_n+1)(Z_n^2+2Z_n+1)~f_3(Z_m)}{(Z_n^2+4Z_n+1)(Z_n^2-1)~f_2(Z_m)}\\
&\simeq \dfrac{3\Delta z}{2}\dfrac{ (Z_n^2+2Z_n+1) ~f_3(Z_m)}{(Z_n^2-1)~f_2(Z_m)}\\
&\simeq \dfrac{3\Delta z}{2}\dfrac{ (Z_n+1) ~f_3(Z_m)}{(Z_n-1)~f_2(Z_m)}
\end{align*}
It can be seen that the oscillatory pole arising out of $Z_n=-1$ has been canceled by one of the repeated zeros introduced in the numerator polynomial by the proposed scheme.

It may be worth recalling here that, invariably the stability analysis of the Galerkin scheme was carried out analytically only for the 1D version of the problem \cite{cdbook}, \cite{up1},  \cite{quadp1}. The associated analysis required solution of the difference equation, which in turn demanded a structured grid in 1D. With respect to the problem in hand, it was the variation of the input perpendicular to the flow was an issue and hence a 2D analysis was required in the above. In order to simplify the analysis, however, an infinite conductor was considered. Practicality requires that the 2D domain must be finite in both $y$ and $z$ directions. Therefore, it becomes necessary to verify the applicability of the analytical findings for the practical situation. 

	\begin{figure}
		\centering
		{\subfigure[]{\label{2dprob2} 
\includegraphics[scale=0.3]{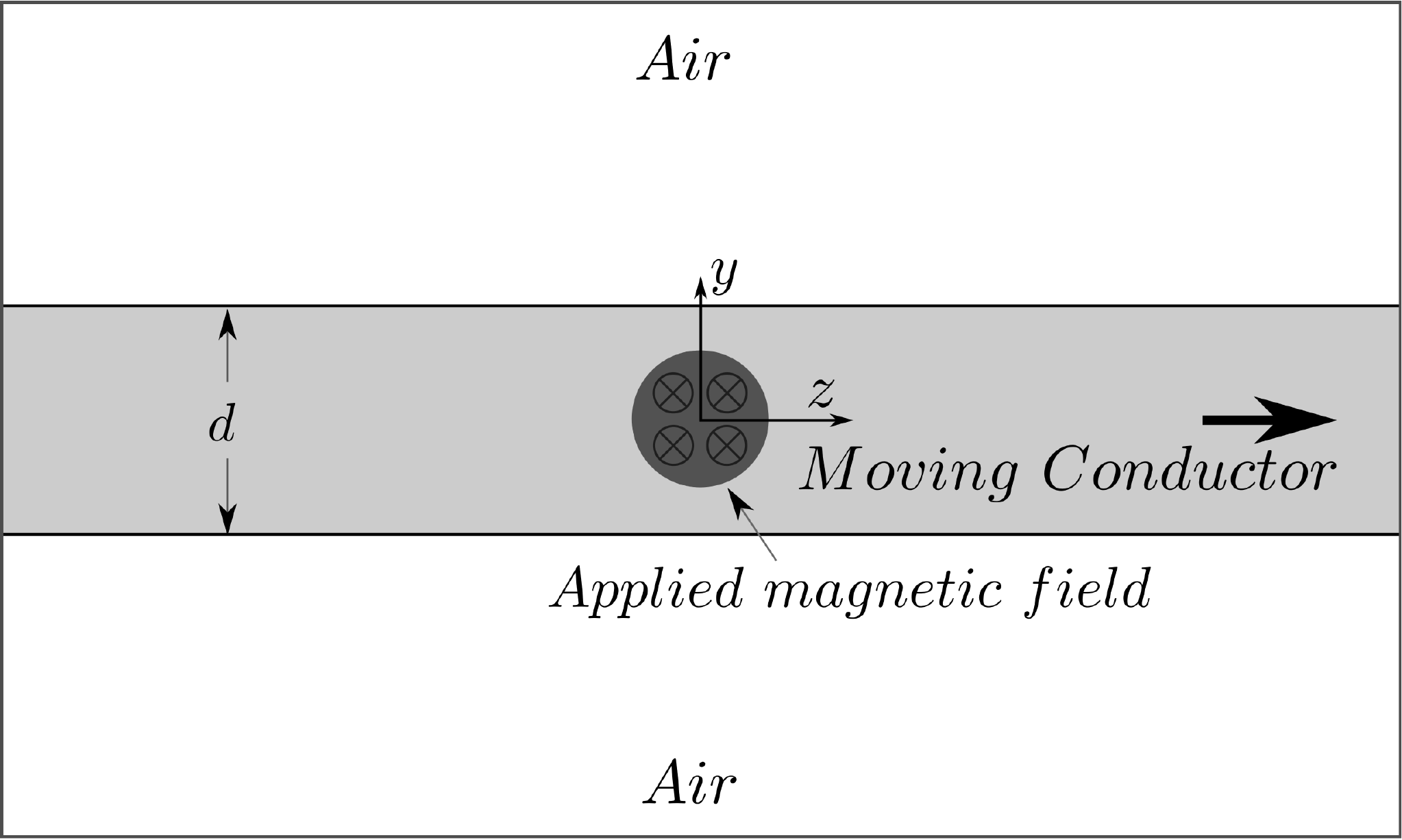}}}		
{\subfigure[]{\label{2dmesh} 
\includegraphics[scale=0.18]{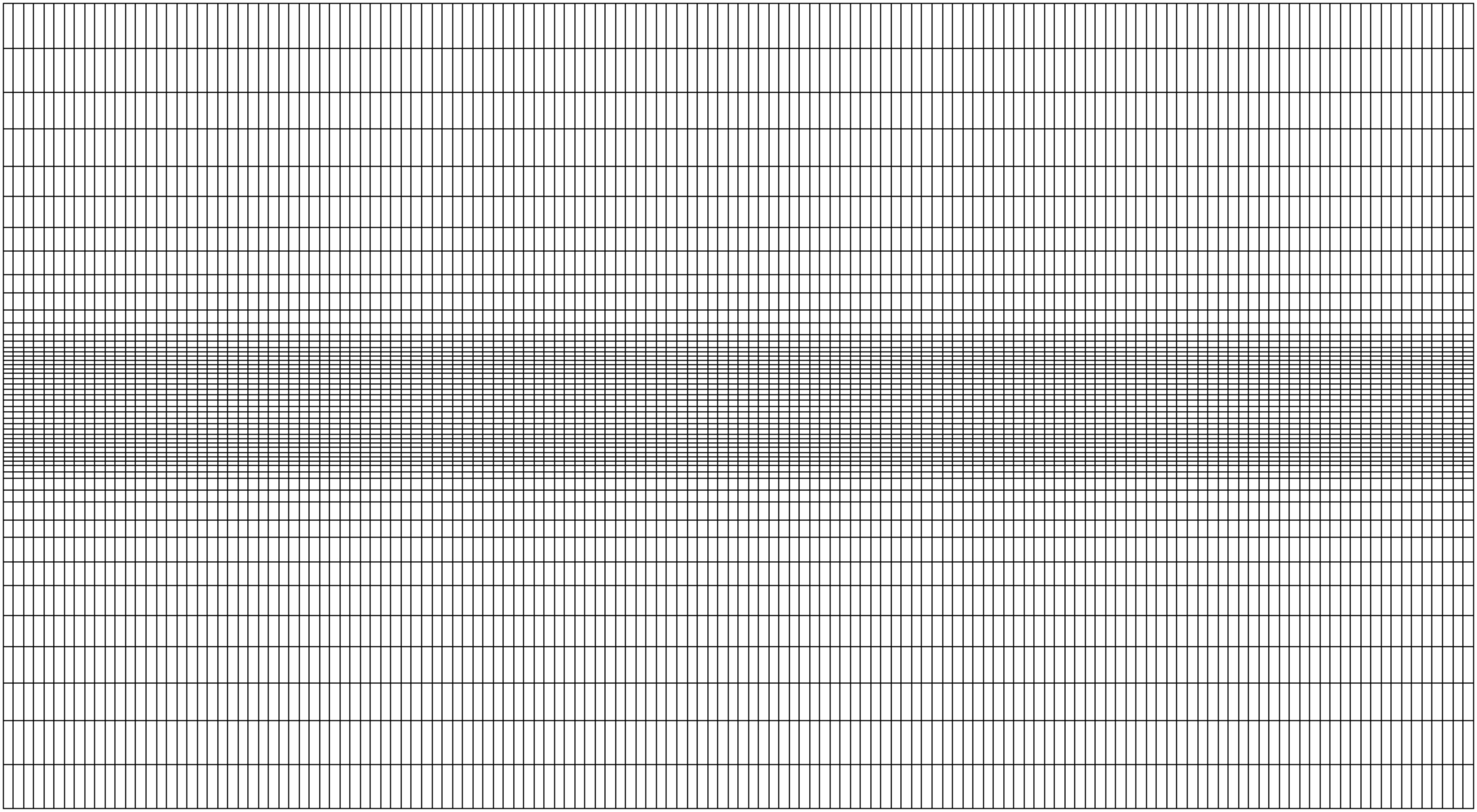}}}
		\caption{2D FEM simulation  $~$  (a) Schematic of the 2D problem (b) 2D mesh}
		\label{un2dresults}
	\end{figure}
For the numerical simulation, problem described in  Fig. \ref{2dprob1} is modified as shown in Fig. \ref{2dprob2}, wherein a conducting sheet of finite thickness ($d$) and conductivity $\sigma=7.21\times10^6Sm^{-1}$ is sandwiched between air regions of thickness $5d$. The axial length is constrained to 6 times the axial width of the input magnetic field. The simulations are carried out with mesh shown in Fig. \ref{2dmesh}, wherein the element size is varied along the $y$-direction. Two different input magnetic field profiles (shown in Figs. \ref{2dipsqr} and \ref{2dipcrc}) are considered.

\begin{figure}
		\centering
		{\subfigure[]{\label{2dsqr1d} 
\includegraphics[scale=0.8]{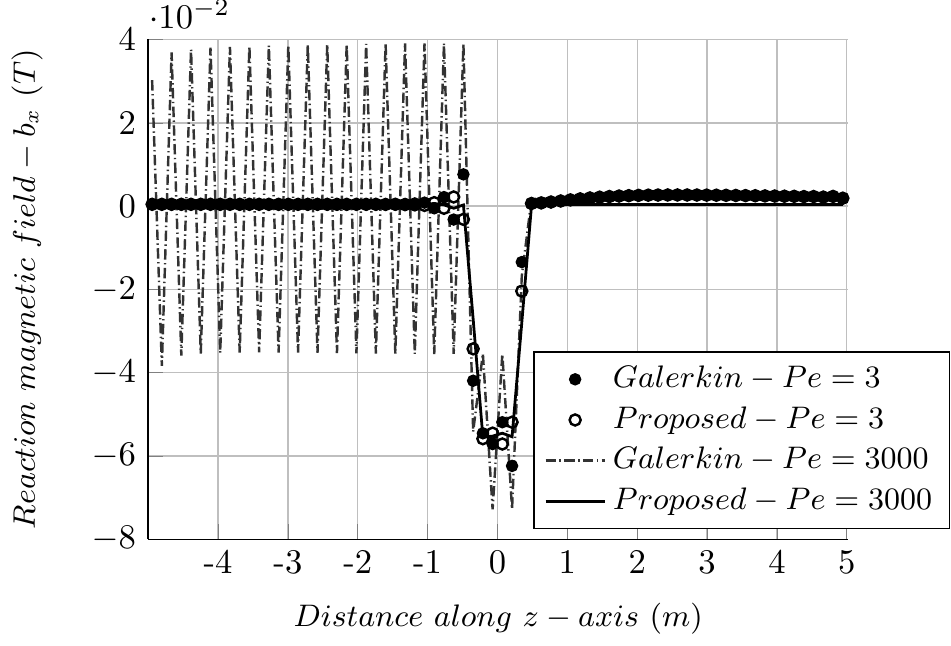}}}		
		{\subfigure[]{\label{2dcrc1d} 
\includegraphics[scale=0.8]{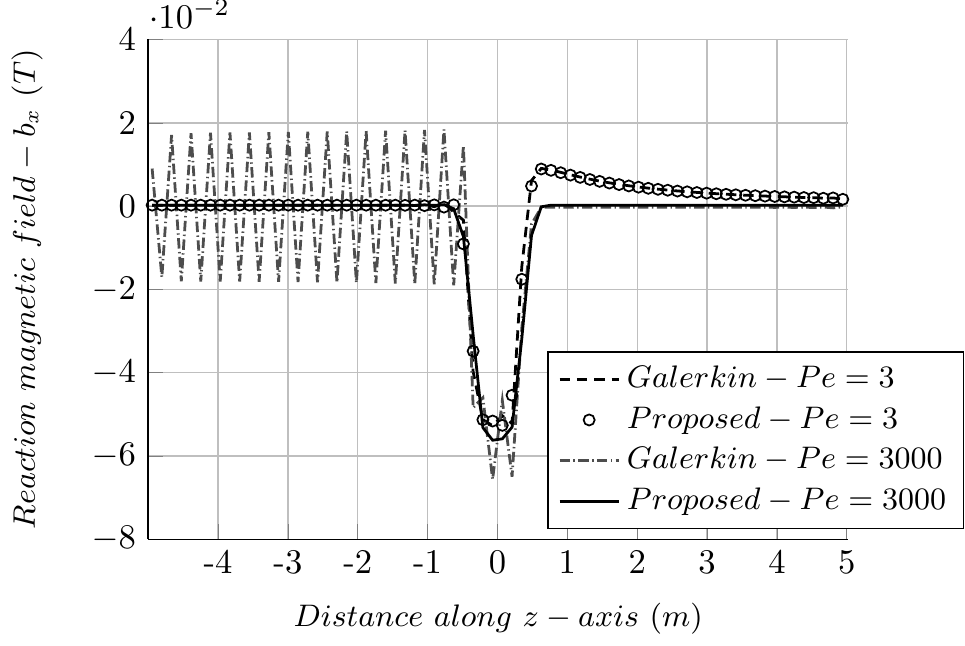}}}
		{\subfigure[]{\label{2dpe60glrk} \includegraphics[scale=0.503]{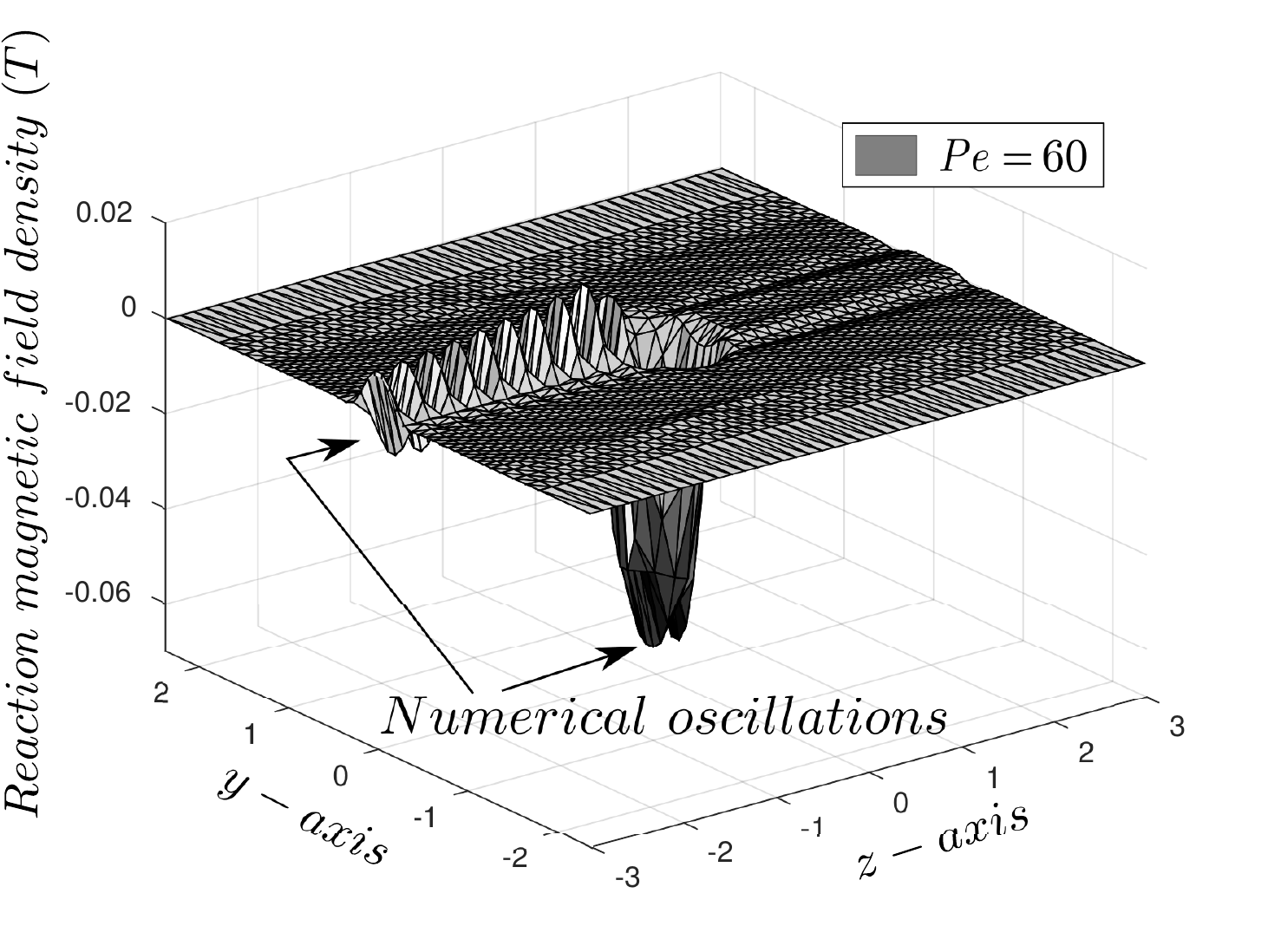}}}
		{\subfigure[]{\label{2dpe60avg} \includegraphics[scale=0.503]{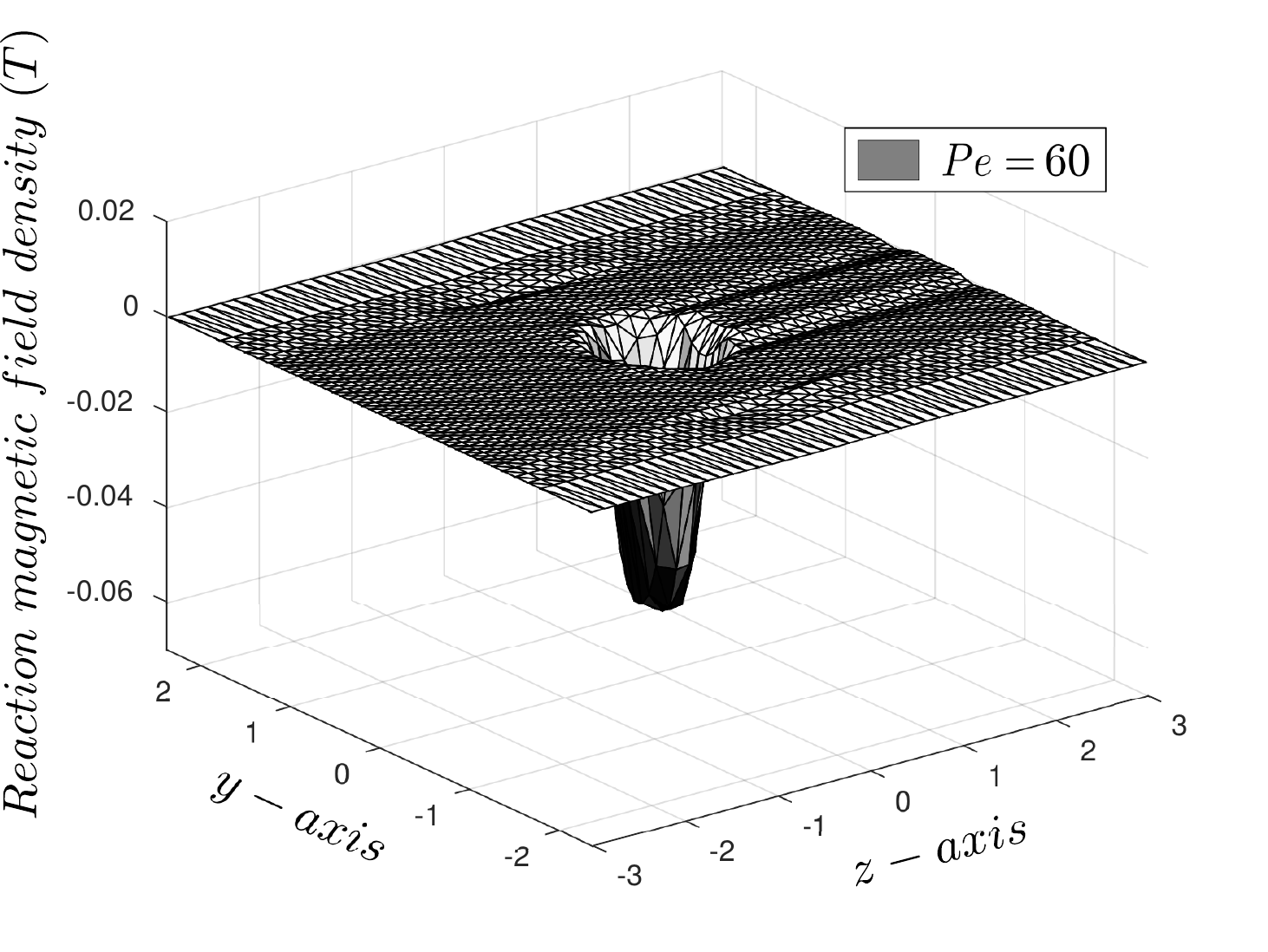}}}
		\caption{Selected results from the 2D FEM simulation  $~$  (a) $b_x$ along the $z$-axis for different Peclet numbers with  input magnetic field of rectangular pulse profile. (b) $b_x$ along the $z$-axis for different Peclet numbers with input magnetic field of smooth circular profile.  (c) Galerkin scheme  with input magnetic field of smooth circular profile.  (d) Proposed scheme with input magnetic field of smooth circular profile. }
		\label{2dresults}
	\end{figure}

Extensive simulations varying almost all the parameters have been carried out. Sample simulation results for $d=1.3m$ and a few selected velocities are presented in Fig. \ref{2dresults}. The computed magnetic field along the axis are presented in Figs. \ref{2dsqr1d} and \ref{2dcrc1d}, for the two different input field profiles and velocities. For further illustration, Figs. \ref{2dpe60glrk} and \ref{2dpe60avg}, present the spatial profile of the field for $Pe=60$. It is evident from the sample results that the proposed scheme, in line with the theoretical prediction, is very stable for $Pe>>1$. However, similar to 1D case, small oscillations prevail in the midrange of $Pe$  ($1<Pe<10$) and it asymptotically vanishes with the increase in $Pe$. 
 
 In order to scrutinize the efficacy of the scheme, the `Testing Electromagnetic Analysis Methods' (TEAM) problem No. 9 \cite{team9} is chosen next. In this problem, the applied magnetic field has components even in the flow direction and further induction to magnetic media is considered.


\subsection{ TEAM problem No. 9}
This problem involves an infinite ferromagnetic material with $\sigma=5\times 10^6 Sm^{-1}$ and $\mu_r=1,50$. This material has a cylindrical bore of diameter $28~mm$. A concentric  current loop of diameter $24~mm$ carrying $1A$  moves at uniform velocity in the bore. This axisymmetric problem is non-uniformly discretised with higher mesh density around the current loop. The FEM model involves 2288 linear quadrilateral elements. For the analysis, worst case involving velocity of $v=100ms^{-1}$ is considered. The resulting $Pe$, due to non uniform discretisation, varies from $5$ to $200$.
	\begin{figure}
		\centering
	  {\subfigure[]{\label{oscbr1} \includegraphics[scale=0.48]{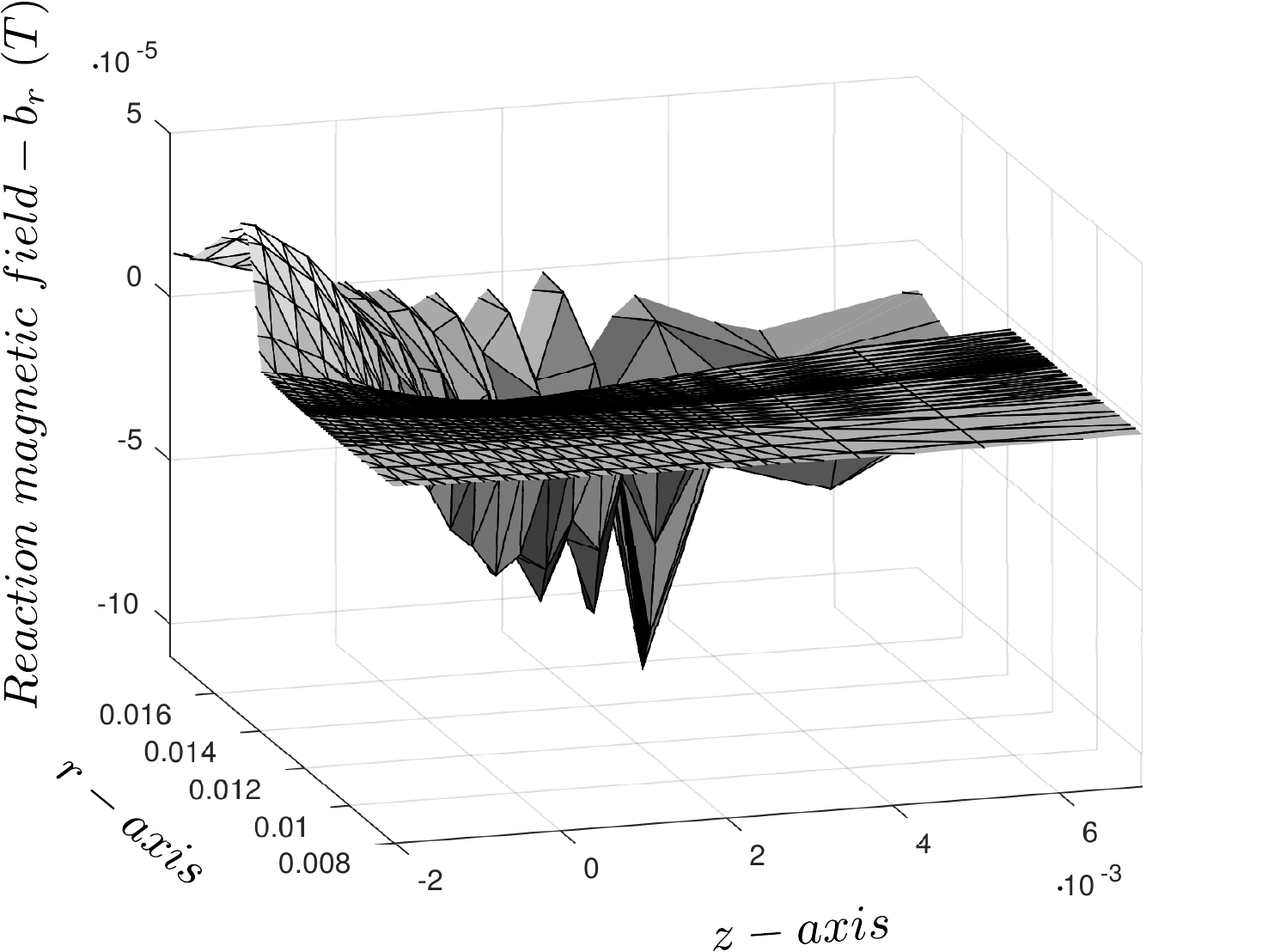}}}
		  {\subfigure[]{\label{stbbr1} \includegraphics[scale=0.48]{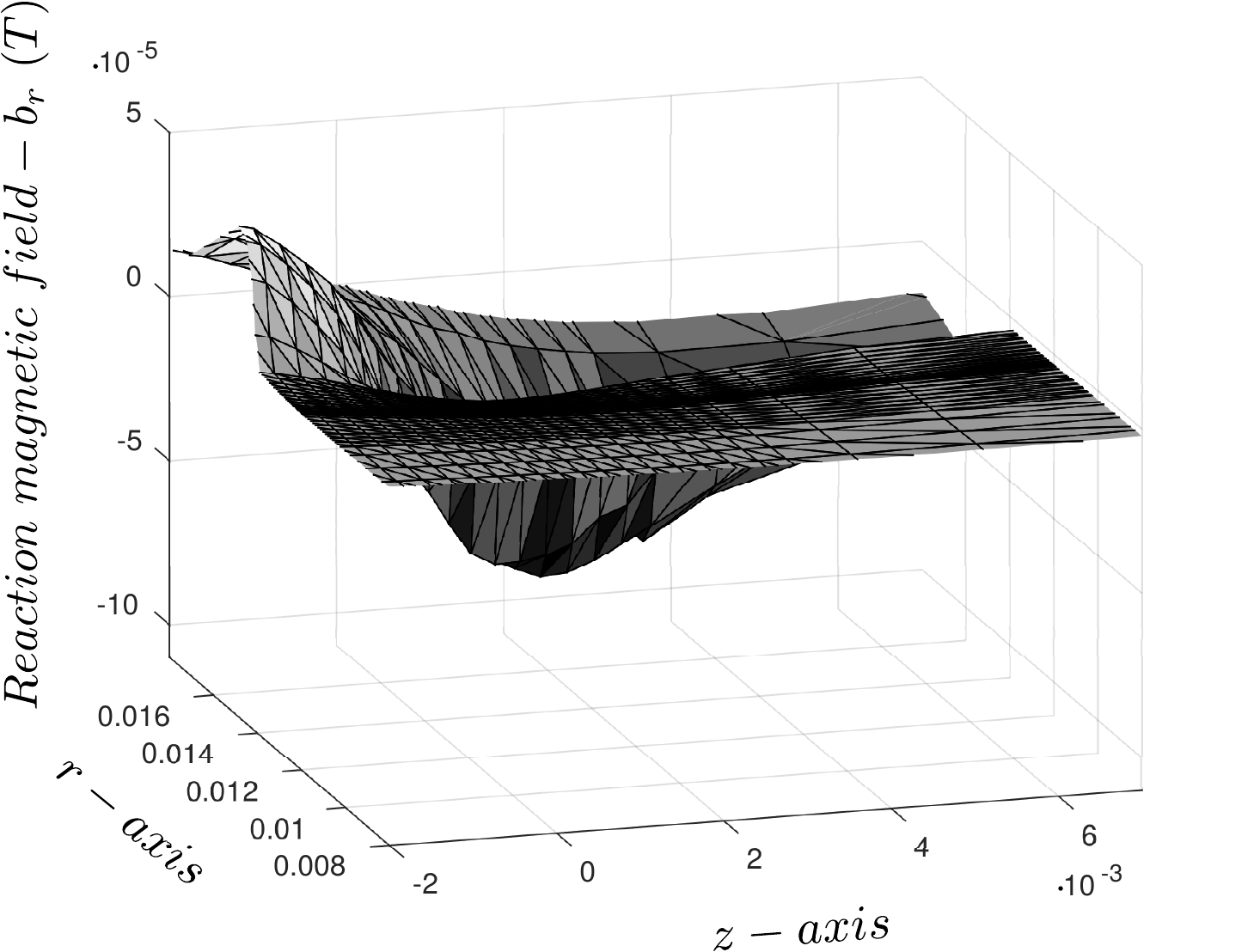}}}
		  {\subfigure[]{\label{team1dair} \includegraphics[scale=0.85]{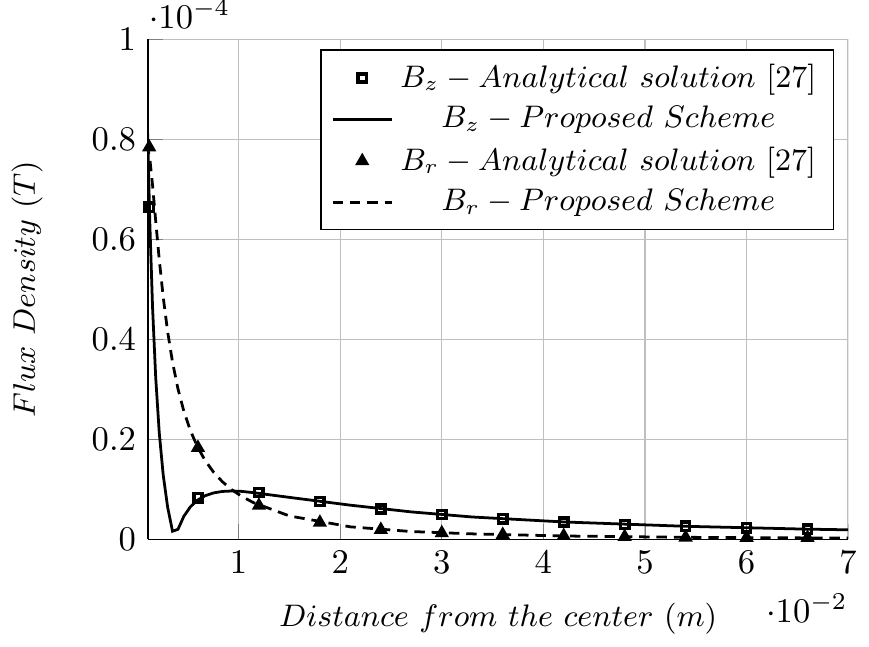} }}
		  {\subfigure[]{\label{team1diron} \includegraphics[scale=0.85]{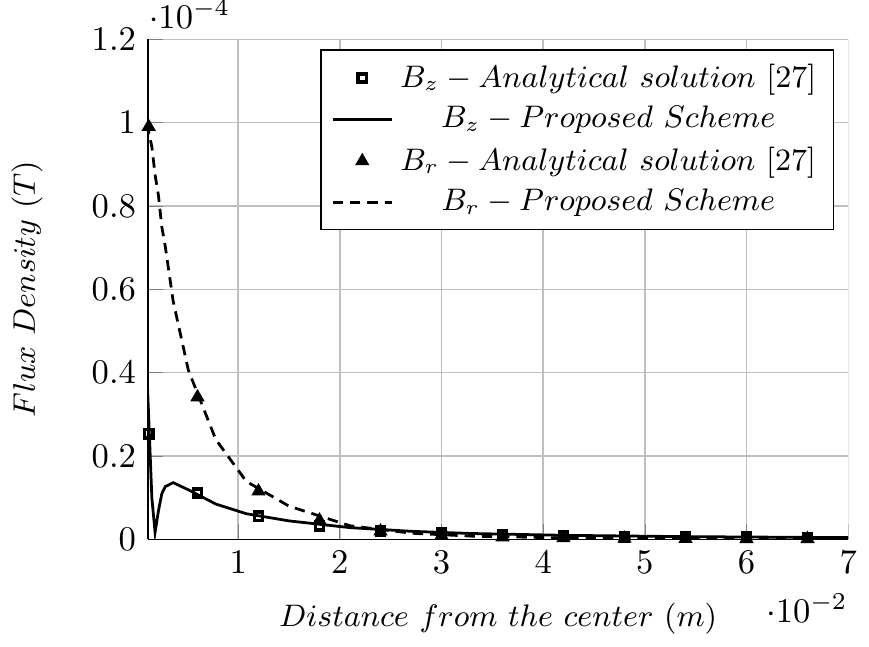}  }}
		\caption{Comparison with TEAM 9 problem (a) Galerkin Scheme - Reaction magnetic field $-~b_r$ for $u_z=100ms^{-1}$ and $\mu_r=50$. (b) Proposed Scheme - Reaction magnetic field $-~b_r$ for $u_z=100ms^{-1}$ and $\mu_r=50$. (c) Comparison of total flux density, $u_z=100ms^{-1}$, $\mu_r=1$.   (d) Comparison of total flux density, $u_z=100ms^{-1}$, $\mu_r=50$. }
		\label{stb3d}
	\end{figure}

It can be seen from the computed magnetic field presented in Figs. \ref{oscbr1} and \ref{stbbr1} that the GFEM leads to oscillations, while the proposed scheme is free of such errors. 
For the quantitative assessment of the accuracy, the analytical results for the radial and the axial air-gap  flux densities presented for $r=13mm$ in \cite{team9} are considered.
From the comparison made in Figs. \ref{team1dair} and \ref{team1diron}, it is evident that the proposed scheme gives accurate results.
It may be cautioned here that even though GFEM gives oscillatory results in the iron region, the air-gap flux densities are not seriously affected and it is also found to give quite accurate results (however is not presented here). 

\subsection{{Simulation for electromagnetic flowmeter}}
\begin{figure}
		\centering
	  {\subfigure[]{\label{magasm} \includegraphics[scale=0.38]{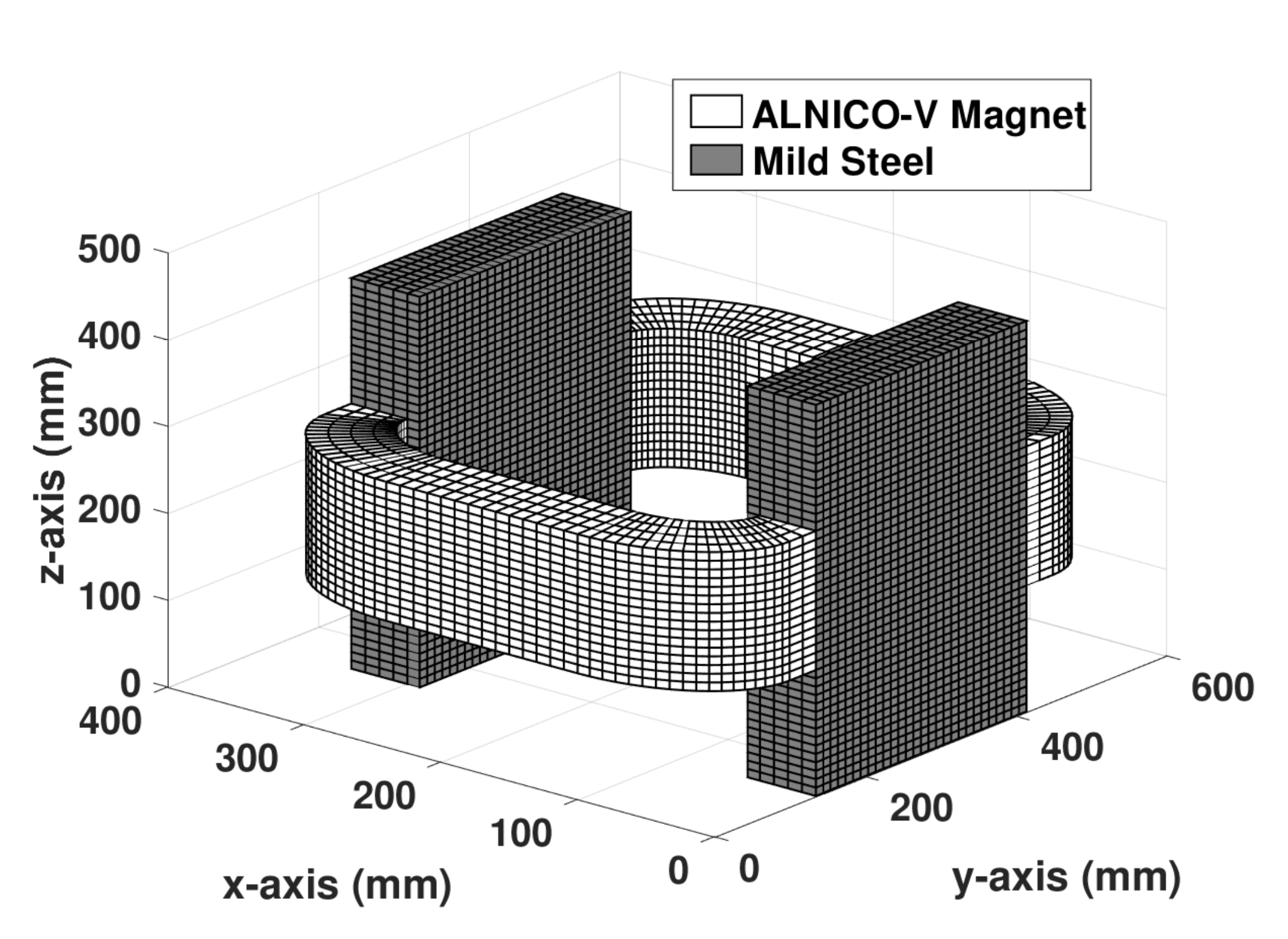}}}
	  	  {\subfigure[]{\label{3ddisc} \includegraphics[scale=0.34]{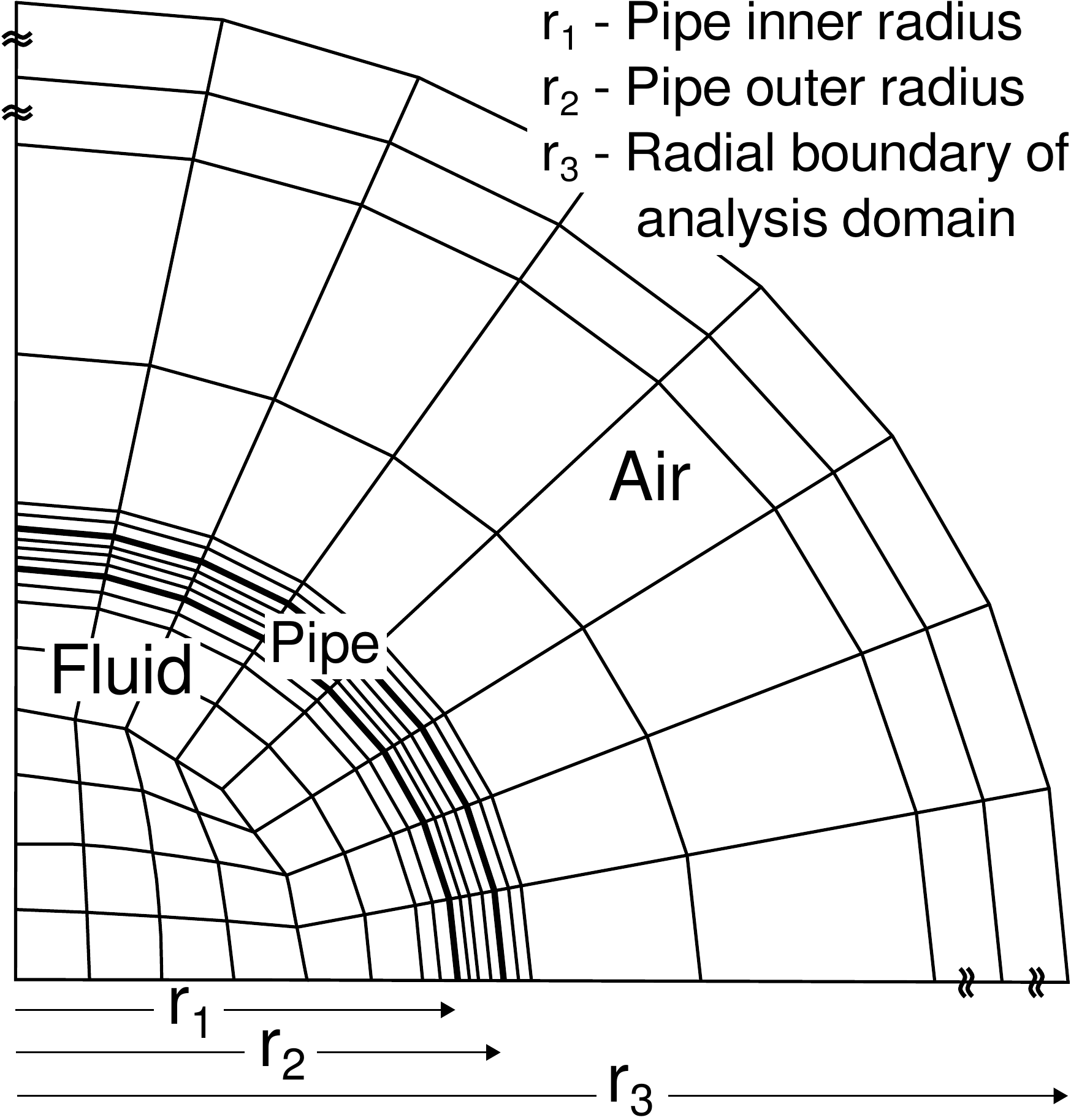}}}
		\caption{{Pertinent details of the flowmeter problem.
 (a) Permanent magnet assembly. (b) Schematic of the discretisation shown here for first quadrant. }}
		\label{3dip}
	\end{figure}

{In order to demonstrate the ability of the proposed scheme in tackling the real life problems, the electromagnetic flowmeter, which is a 3D problem, is simulated for high Peclet numbers. The input magnetic field is provided by the permanent magnet assembly shown in Fig. {\ref{magasm}}. The ambient magnetic field is separately evaluated using fictitious magnetic charge method} {\cite{ukreport}} {and its peak value is about $0.0413 ~T$. The outer diameter of the steel pipe is $0.2191~ m$ and it has a thickness of $0.0164~m$. Liquid sodium at temperature $200^0 C$ is considered as the flowing fluid.}

{Discretisation of the 3D geometry is carried out with 8 node brick elements spanning radially into the air region surrounding the pipe. Axially, analysis domain extends up to $9.5~m$, which corresponds to 21 times the length of the magnetic assembly. Totally $31211$ nodes and $29696$ elements are employed. A schematic of the discretisation in the cross section is given in Fig. {\ref{3ddisc}} and structured meshing is employed along the flow direction.}

\begin{figure}
		\centering

		  {\subfigure[]{\label{caseA} \includegraphics[scale=1]{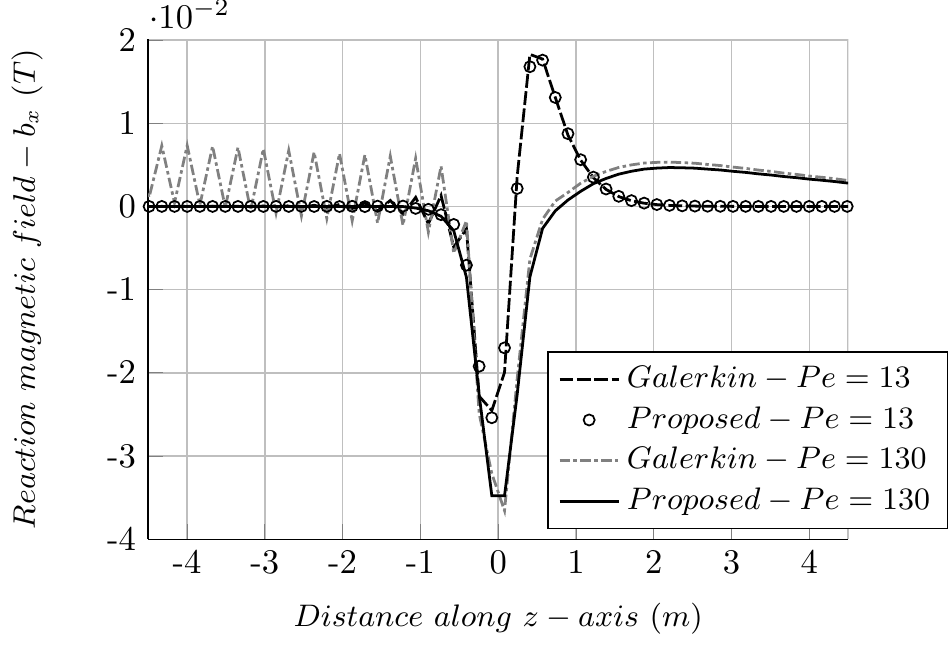}}}
		  {\subfigure[]{\label{glrk100nb} \includegraphics[scale=0.143]{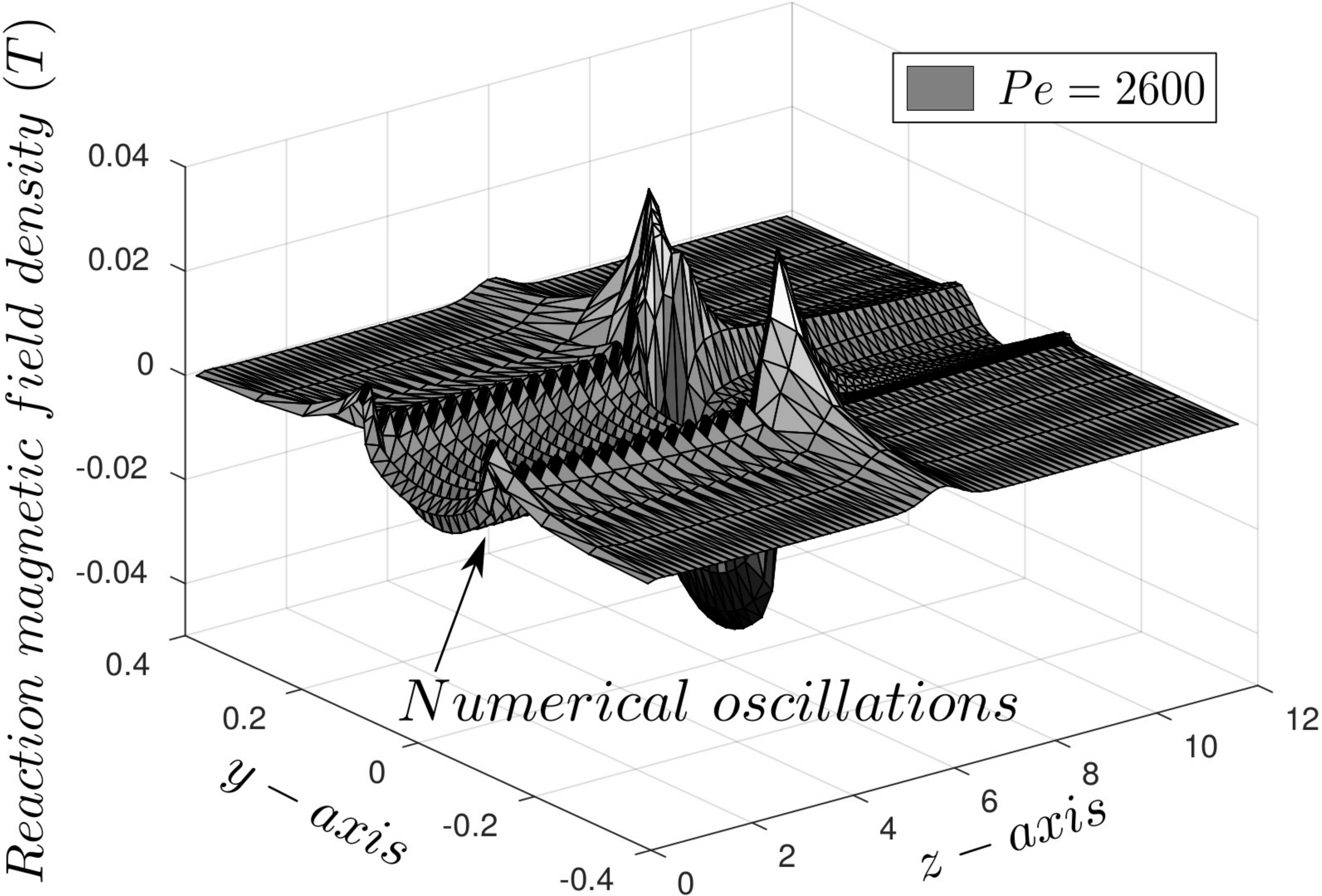} }}
		  {\subfigure[]{\label{avg100nb} \includegraphics[scale=0.143]{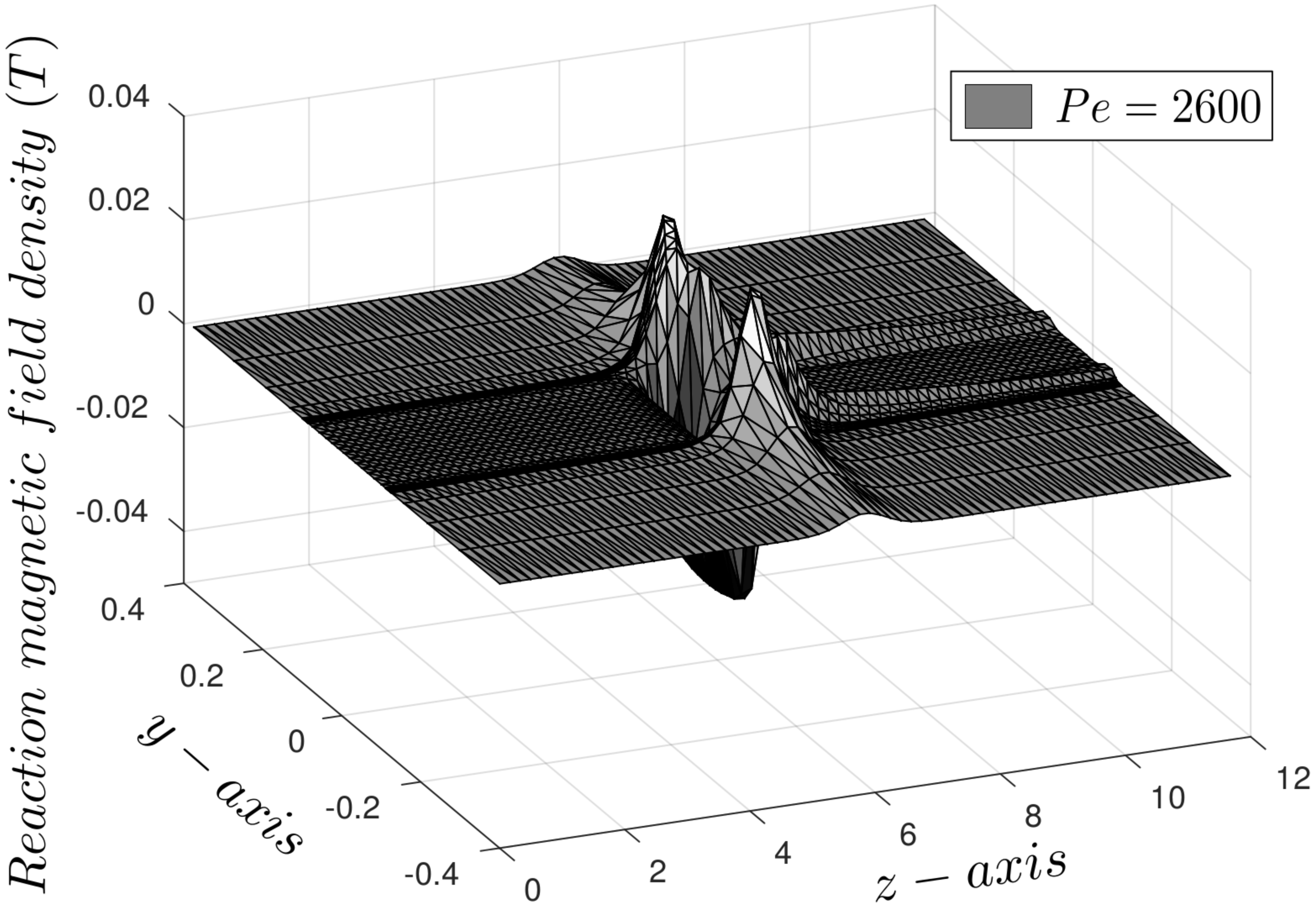}  }}
		  		\caption{{Sample simulation results for electromagnetic flowmeter (a)  $x$-component of ${\bf{b_{r}}}$ along the pipe axis ($z$-axis) for different Peclet numbers (b) Galerkin Scheme - $x$-component of reaction magnetic field ${\bf{b_{r}}}$ in  $~x$ = 0 plane  $~$ (c)  Proposed Scheme - $x$-component of reaction magnetic field ${\bf{b_{r}}}$ in  $~x$ = 0 plane.}}
		\label{3dip}
	\end{figure}

{As per the proposed scheme, flux density averaged over the eight nodes defining the element, was employed for the evaluation of the elemental matrices (refer to  ({\ref{3dipbxyz}})).}
\begin{align}
B_{xe} = \dfrac{1}{8} \sum\limits_{n=1}^{8} B_{xn} ;  ~~~B_{ye} = \dfrac{1}{8} \sum\limits_{n=1}^{8} B_{yn} ; ~~~B_{ze} = \dfrac{1}{8} \sum\limits_{n=1}^{8} B_{zn} ;
\label{3dipbxyz}
\end{align}
{Simulations are carried out for a wide range of $Pe$ and sample results comparing the present scheme with the GFEM is provided in Figs.
 {\ref{caseA}}, {\ref{glrk100nb}} and {\ref{avg100nb}}. It is evident from the figures that the proposed scheme provides a stable and accurate results even for  the flow rates well beyond the practical operating range. }

\subsection{Discussion}

{Simulations were carried out for several problems, including 3D cases and only the representative ones were presented in the earlier sections. In all these cases, accuracy as determined by taking GFEM results with very fine discretisation as the reference, is found to be quite high.

But for the small modification to the RHS term, the procedure adopted is identical to GFEM. As a result, the increase in computational time with the present approach is observed to be less than 2\% for 2D and 3D cases. }

{It may be recalled that, in the definition of the problem itself, a kind of graded regular mesh involving elements of equal lengths in the flow direction was envisaged for every layer. However, layer to layer length could be different.
In order to make a quick assessment of the arbitrariness in the general FEM mesh, two kinds of numerical experiments are conducted with a 2D problem shown in Fig. {\ref{2dprob2}}. In the first, skewness is introduced for every element in the graded regular mesh and in the second, arbitrary mesh is employed. In the former, the loss of accuracy was less than 5 \%, even when the interior angle was changed to  $90\pm45^0$. In the case of the latter, accuracy was marginally affected at lower $Pe$ ($< 100$) and significant deterioration was found thereafter.}

Incidentally, successful efforts have also been  made to extend the philosophy to 1D second order elements. It is then intuitively extended to second order 2D and 3D elements. However, due to page restrictions, they are not dealt here.

%
%

\section{Summary and Conclusion}
The GFEM, when employed for magnetic problems involving conductor moving at high velocities is known to suffer from numerical instability.  To address this problem, SU/PG scheme is generally suggested. However the SU/PG scheme is known to suffer from error at the boundary transverse to the velocity. Recently a GFEM based, simple, stable scheme has been proposed. However, in that the input field is assumed to vary only along the direction of the velocity
 \cite{su1}. To overcome this limitation a simple alternative approach is proposed in this work.

Similar to our earlier work, the problem is analyzed in the Z-transform domain and the pole-zero cancellation principle is adopted to propose a new stable scheme. This involved, restatement of the input flux density as an elemental weighted average. Analytically, it is shown for both 1D and 2D versions of the problem, that the proposed scheme is absolutely stable at high values of $Pe$.
Incidentally, at the mid-range of $Pe$, the proposed scheme exhibits small oscillations, the upper bound for which is theoretically shown to be  $<5\%$.
These predictions are adequately validated with numerical exercises, which included No. 9 of the TEAM benchmark problem.

\bibliographystyle{IET} 
\bibliography{References}

\section{Appendix}
\subsection{{Details of the Z Transforms}}
{The literature, including text books are abundant on Z-transform analysis. In order to provide a quick reference, especially for 2D Z-transform, some of the relevant aspects are reproduced here from }\cite{dcbook1, dcbook2, 2dhb, 2dsi, 2dsb}.

\subsubsection{{1D Z transform definition}}
{For an equi-spaced discrete sequence {$X_{[n]}$}, the Z-transform  is defined as }{\cite{dcbook1, dcbook2, 2dhb}},
\begin{equation}
X(Z) = Z[X_{[n]}] = \sum\limits_{n=-\infty}^\infty X_{[n]} Z^{-n}
\end{equation}
{The Z-transform possess the shifting property given by }\cite{dcbook2},
\begin{equation}
Z[X_{[n \pm k]}] =Z^{ \pm k} X(Z)
\end{equation}
{The RHS of the linear difference equation can be related to the variable in the LHS, using transfer function in Z-domain} \cite{dcbook2}. { The transfer function in the Z-domain {$H(Z)$} is defined as,}
\begin{equation}\label{hztf}
H(Z) = \dfrac{U(Z)}{V(Z)} = b_0 \dfrac{\prod\limits_{s=1}^{S} (Z-b_s)}{\prod\limits_{t=1}^{T} (Z-a_t)},   ~~~~ s \leq t
\end{equation}
{For the system to be stable, the poles of the transfer function ({\ref{hztf}}) in  Z-plane must lie within the unit circle ({$|a_t|   ~\textless~ 1$}) which defines the region of stability} \cite{dcbook1}.
\subsubsection{{2D Z transform definition}}
{For a 2D discrete sequence {$X_{[n,m]}$}, the 2D Z-transform is defined as} \cite{2dhb},
\begin{equation}
X(Z_n,Z_m) = Z[X_{[n,m]}] = \sum\limits_{n=-\infty}^\infty \sum\limits_{m=-\infty}^\infty X_{[n,m]} Z^{-n}Z^{-m}
\end{equation}
{It also exhibits shifting property given by} \cite{2dhb},
\begin{equation}
Z[X_{[n \pm k,m \pm l]}] =Z_n^{\pm k} Z_m^{\pm l} X(Z_n,Z_m)
\end{equation}
{ Similar to 1D case,  transfer function {$H(Z_n,Z_m)$}, can be defined for the 2D case} {\cite{2dsi}}. { Under special circumstances, like the one encountered in the present work, it assumes a separable form as given below } {\cite{2dsi}} {\cite{2dsb}}
\begin{equation}\label{hztf2d}
H(Z_n,Z_m) = \dfrac{U(Z_n,Z_m)}{V(Z_n,Z_m)} = b_{0n}b_{0m} \dfrac{\prod\limits_{sn=1}^{Sn} (Z_n-b_{sn}) \prod\limits_{sm=1}^{Sm} (Z_m-b_{sm})  }{\prod\limits_{tn=1}^{Tn} (Z_n-a_{tn}) \prod\limits_{tm=1}^{Tm} (Z_n-a_{tm})},   ~~~~ sn \leq tn,  sm \leq tm 
\end{equation}
{For such transfer functions, the stability regions are unit circle in their respective Z-planes and hence the poles must lie within them ({$|a_{tn}|   ~\textless~ 1$} and {$|a_{tm}|   ~\textless~ 1$}) .} {On the other hand, if the pole lies  on the circumference only marginal stability is ensured. Further, if the pole is at -1, then $U(Z_n,Z_m)$ will exhibit sustained oscillation} {\cite{dcbook1, dcbook2}.
\subsection{Analytical solution of the difference equation}
	\begin{figure}
		\centering
		{\subfigure[]{\label{gal_diff:fig} \includegraphics[scale=0.45]{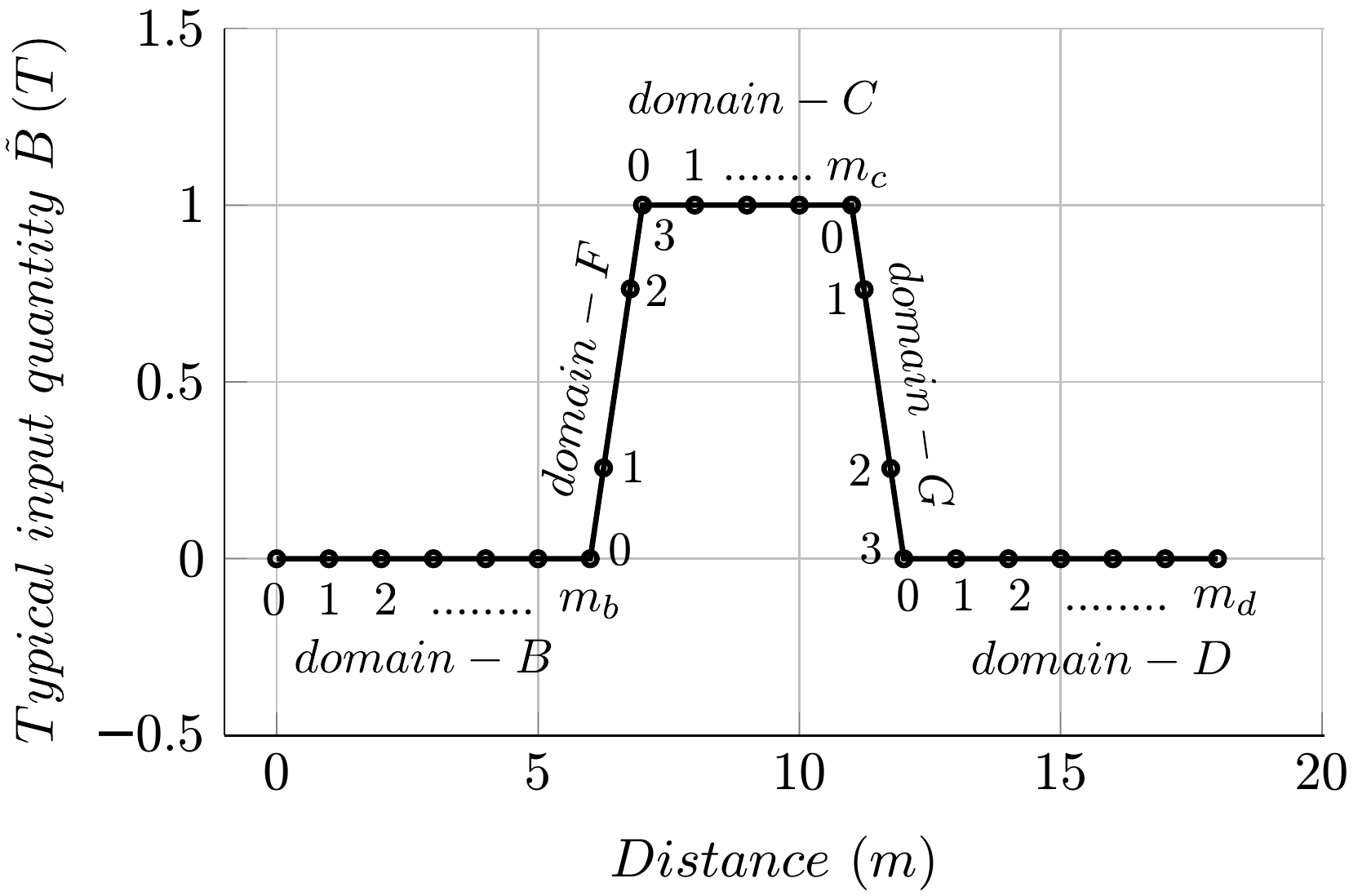}}}
		{\subfigure[]{\label{gal_diff:g1} \includegraphics[scale=0.8]{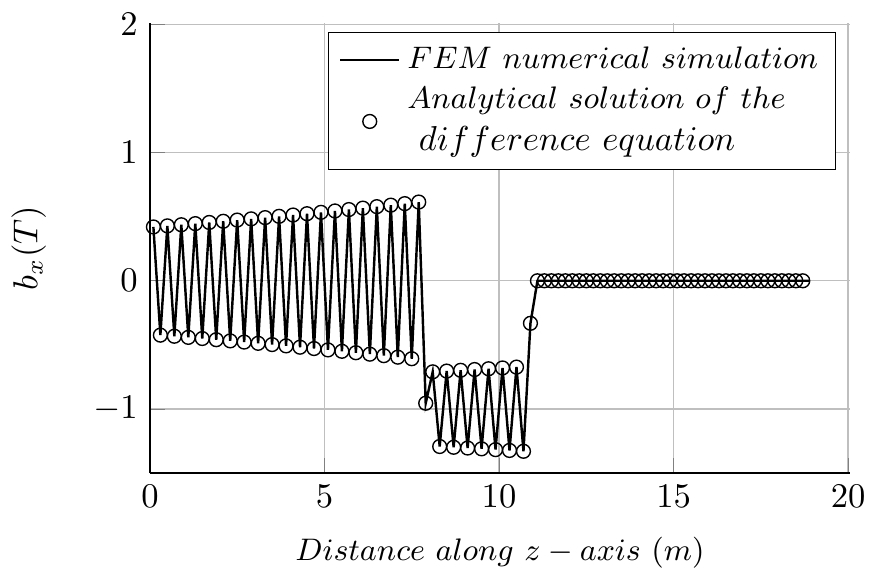}}}
		  {\subfigure[]{\label{gal_diff:p1} \includegraphics[scale=0.8]{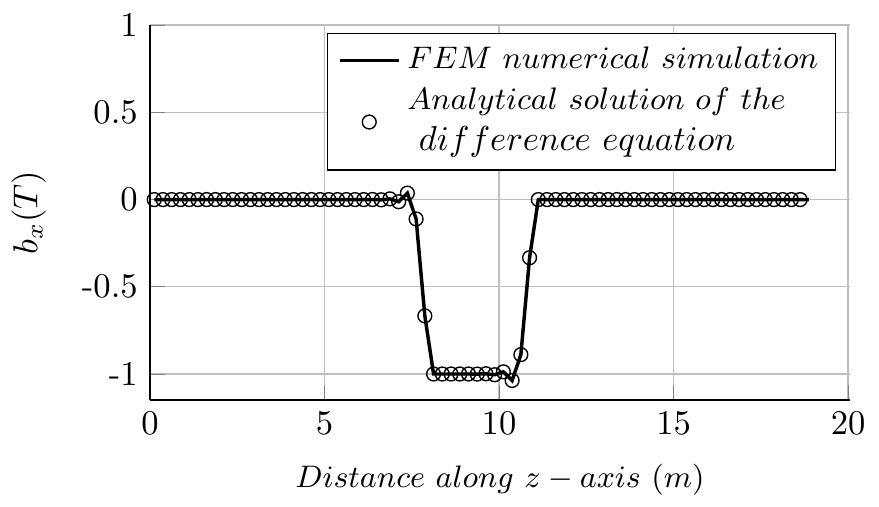}}}
		  {\subfigure[]{\label{gal_diff:p2} \includegraphics[scale=0.8]{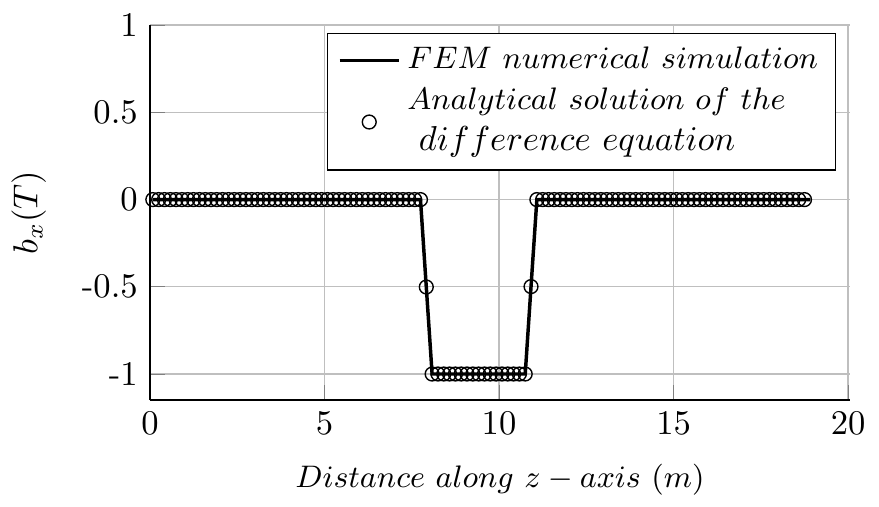}}}
		\caption{Analytical solution of the difference equation (a) Five sub-domains and their ranges (b) Validation of the analytical solution for the \emph{Galerkin scheme}  $~$ $Pe=200,~\Delta z = 0.20,~ m_c=12, ~m_b = 38, ~m_d = 38.~~$ (c) Validation of the analytical solution for the \emph{proposed scheme} $~$ $Pe=2,~\Delta z = 0.25,~ m_c=9, ~m_b = 30, ~m_d = 30 $.  (d) Validation of the analytical solution for the \emph{proposed scheme} $~$ $Pe=400,~\Delta z = 0.17,~ m_c=15, ~m_b = 46, ~m_d = 46 $.}
		\label{gal_diff}
	\end{figure}
The solution domain is divided into five sub-domains ($B,~C,~D,~F,~G$) as shown in Fig. \ref{gal_diff:fig}  where $\tilde{B} = (B_{n-1}+4B_n+B_{n+1})/6$ for Galerkin scheme and $\tilde{B} = (B_{n-1}+2B_n+B_{n+1})/4$ for the  proposed scheme. The solutions of the different domains are \cite{su1},

Domain $B$ ($0 \leq n_b \leq m_b $): $~~y_b(n_b)=b_1+b_2r^{n_b}  $

Domain $F$ ($0 \leq n_f \leq 3 ~~$): $~~y_f(n_f)=f_1+f_2r^{n_f}+y_{pf}(n_f)$

Domain $C$ ($0 \leq n_c \leq m_c$): $~~y_c(n_c)=c_1+c_2r^{n_c}+y_{pc}(n_c)$

Domain $G$ ($0 \leq n_g \leq 3~~$): $~~y_g(n_g)=g_1+g_2r^{n_g}+y_{pg}(n_g)$

Domain $D$ ($0 \leq n_d \leq m_d$): $~~ y_d(n_d)=d_1+d_2r^{n_d}$

where, $r = (-1-Pe)/(-1+Pe), ~y = A_y, y_p -$particular solution,   suffixes $b,c,d,f,g$ designate the domain names and $b_1$, $b_2$,  $c_1$, $c_2$,  etc. are the parameters of the complimentary solution belonging to their respective domains. 
By imposing the boundary conditions, two parameters are found.
\[y_b(0)=0 ~~\Rightarrow b_1=-b_2; ~and ~ y_d(m_d) = y_d(m_d+1)~~~\Rightarrow~~ d_2 = 0\]
Then, by imposing the equality condition and by satisfying the difference equation at the joining nodes, the other parameters of the complementary solutions are found \cite{su1}.
\begin{flushleft}
\begin{align}
g_2 &= \dfrac{r(y_{pg}(2)-y_{pg}(3))}{r^{m_g}(r-1)} \label{g2}\\
c_2 &= \frac{r(y_{pc}(m_c-1)-y_{pc}(m_c))}{r^{m_c}(r-1)} + \dfrac{g_2}{r^{m_c}}+\dfrac{y_{pg}(1)}{r^{m_c}(r-1)}+\dfrac{\lambda}{r^{m_c}}\label{c2}\\
f_2 &=\frac{r(y_{pf}(2)-y_{pf}(3))}{r^{m_f}(r-1)} + \dfrac{c_2}{r^{m_f}}+\dfrac{y_{pc}(1)}{r^{m_f}(r-1)}+\dfrac{\lambda}{r^{m_f}}\label{f2}\\
b_2 &= \dfrac{f_2}{r^{m_b}}+\dfrac{y_{pf}(1)}{r^{m_b}(r-1)}\label{b2}\\
f_1 &= b_1 + b_2 r^{m_b} - f_2\label{f1}\\
c_1 &= f_1 + f_2 r^3 + y_{pf}(3) -c_2\label{c1}\\
g_1 &= c_1 + c_2 r^{m_c} + y_{pc}(m_c) - g_2\label{g1}\\
d_1 &= g_1 + g_2 r^3 + y_{pg}(3) \label{d1}
\end{align}
\end{flushleft}
The particular solutions for domain $C$ is same for both the schemes,
\[ y_{pc}(n_c) = \lambda n_c \] where, $\lambda = B \Delta z$.
However, the particular solutions in domains $F~\&~G$ are different.

For the Galerkin scheme,
\[y_{pf}(n_f)= -\dfrac{\lambda}{24}n_f^4 + \dfrac{\lambda(r-2)}{6(r-1)}n_f^3 + \dfrac{\lambda(r^2 - 5)}{8(r-1)^2}n_f^2 -  \dfrac{\lambda(r^3 - 10r^2 + 17r + 4)}{12(r-1)^3}n_f\]
\[y_{pg}(n_g)= \dfrac{\lambda}{24}n_g^4 - \dfrac{\lambda(r-2)}{6(r-1)}n_g^3 - \dfrac{\lambda(r^2 - 5)}{8(r-1)^2}n_g^2 +  \dfrac{\lambda(13r^3 - 46r^2 + 53r - 8)}{12(r-1)^3}n_g\]

and for the proposed scheme,
\[y_{pf}(n_f)= -\dfrac{\lambda}{48}n_f^4 + \dfrac{\lambda(r-2)}{12(r-1)}n_f^3 + \dfrac{\lambda(7r^2 -8r- 11)}{48(r-1)^2}n_f^2 +  \dfrac{\lambda(r^3 + 8r^2 - 19r -2)}{24(r-1)^3}n_f\]
\[y_{pg}(n_g)= \dfrac{\lambda}{48}n_g^4 - \dfrac{\lambda(r-2)}{12(r-1)}n_g^3 - \dfrac{\lambda(7r^2 -8r- 11)}{48(r-1)^2}n_g^2 +  \dfrac{\lambda(23r^3 - 80r^2 + 91r -22)}{24(r-1)^3}n_g\]

As a validation, the analytical solutions of the difference equations are compared with the numerical solutions obtained from the FEM in Figs. $\ref{gal_diff:g1}$ and $\ref{gal_diff:p1}$ and $\ref{gal_diff:p2}$.
\subsubsection{Evaluation of the error}
{In }{\cite{su1}}{ the error in the numerical solution for a different numerical scheme was quantified. Following the same procedure, the error in the present numerical scheme is quantified by comparing it with the analytical solution and it takes the form,}
\begin{equation}\label{oscan2}
\widehat{b} =\frac{c_2(r^{m_c-1}-r^{m_c})}{\Delta z}\\
\end{equation}
 Substituting (\ref{c2}) in (\ref{oscan2}), the error in the proposed scheme is found to be,
\begin{equation}\label{errp}
\widehat{b_{p}} = \dfrac{B(r^2 + 2r + 1)}{4r^3}= \frac{B(1-Pe)}{(1+Pe)^3} \\
\end{equation}	
and for the Galerkin scheme,
\begin{equation}\label{errg}
\widehat{b_{g}} =\dfrac{B(r^2 + 4r + 1)}{6r^3} = \frac{B(Pe^2 - 3)(Pe - 1)}{3(Pe + 1)^3}\\
\end{equation}
The extremum of (\ref{errp}) gives the location and the value of the peak error in the proposed scheme.


\begin{center}
-------------------------------------------------
\end{center}


\vspace{10mm}

\begin{minipage}{0.17\linewidth}
{\includegraphics[width=1in,height=1.25in,clip,keepaspectratio]{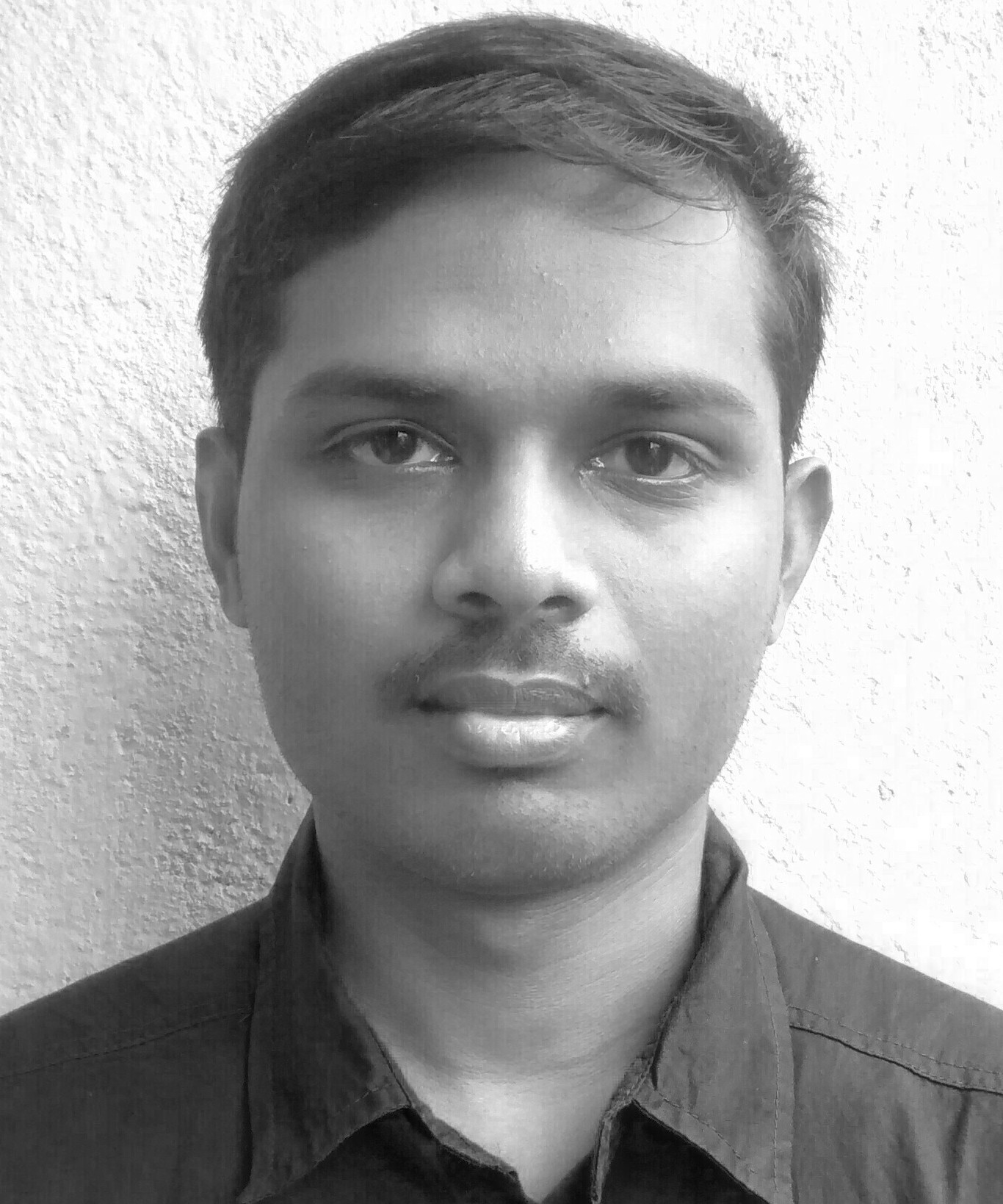}}
\end{minipage}
\begin{minipage}{0.84\linewidth}
{\bf Sethupathy Subramanian}
received the bachelors degree in electrical and electronics engineering from the Anna University, Chennai, India, in 2009. He received the masters degree in electrical engineering from the Indian Institute of Science, Bangalore, India in 2011, where he is currently pursuing his Ph.D. degree. 

$~~~~~~$His research interests include Electromagnetism, Magnetohydrodynamics  and Computational Methods.
\end{minipage}
%
\begin{minipage}{0.18\linewidth}
{\includegraphics[width=1in,height=1.25in,clip,keepaspectratio]{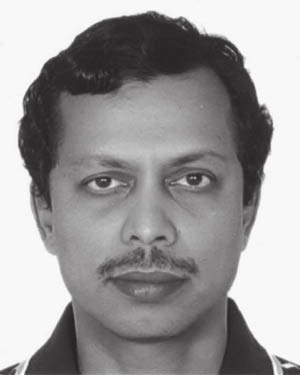}}
\end{minipage}
\begin{minipage}{0.82\linewidth}
{\bf Udaya Kumar}
 was born in Udupi District, Karnataka, India, in 1966. He received the bachelor’s degree in electrical engineering from Bangalore University, Bangalore, India, in 1989 and the M.E. and Ph.D. degrees in high-voltage engineering from the Indian Institute of Science, Bangalore, in 1991 and 1998, respectively. 
\end{minipage}

\vspace{2mm}
He was a Senior Analyst with the Electromagnetic Group of Electromagnetic Research Consultants, Bangalore. Since 1998, he has been with the Indian Institute of Science, where he is currently a Professor in the Department of Electrical Engineering. His research interests include lightning, electromagnetism, and high-frequency response of windings.

He is a member of CIGRE working groups WG C4.26 on “Evaluation of Lightning Shielding Analysis Methods for EHV and UHV DC and AC Transmission Lines” and WG C4.37 on “Electromagnetic Computation Methods for Lightning Surge Studies with Emphasis on the FDTD Method.”


\end{document}